\documentclass[12pt,onecolumn]{IEEEtran} % arxiv Mode
\IEEEoverridecommandlockouts                              % This command is only
\usepackage[top = 1.91cm, bottom = 1.91cm, left = 1.91cm, right = 1.91cm]{geometry}
\usepackage{url}
\usepackage{afterpage}
\usepackage[pdftex]{graphicx}
\usepackage{caption}
\usepackage{amsmath} % assumes amsmath package installed^{}
\usepackage{amssymb}  % assumes amsmath package installed

\usepackage{ntheorem}
\usepackage{eurosym}
\usepackage{color}
\usepackage{ae}
\usepackage{cite}
\usepackage{acronym}
\usepackage[detect-all]{siunitx}

\newcommand{\comment}[1]{}

\theoremstyle{plain}
\theoremheaderfont{\normalfont\bfseries}
\theorembodyfont{\normalfont}
\theoremseparator{:}
\title{{Analyzing Operational Flexibility \\ of Electric Power Systems}}

\author{
	Andreas Ulbig and G\"{o}ran Andersson\\
	Power Systems Laboratory, ETH Zurich, Switzerland\\
	{\textmd ulbig$\;$\textbar$\;$andersson$\,$@$\,$eeh.ee.ethz.ch}
}

\begin{document}
%\newgeometry{top = 2.54cm, bottom = 1.91cm, left = 1.91cm, right = 1.91cm}

\maketitle

\begin{abstract}
Operational flexibility is an important property of electric power systems and plays a crucial role for the transition of today's power systems, many of them based on fossil fuels, towards power systems that can efficiently accommodate high shares of variable \ac{res}. The availability of sufficient operational flexibility in a given power system is a necessary prerequisite for the effective grid integration of large shares of fluctuating power in-feed from variable \ac{res}, especially wind power and \ac{pv}. 
This paper establishes the necessary framework for quantifying and visualizing the technically available operational flexibility of individual power system units and ensembles thereof.
Necessary metrics for defining power system operational flexibility, namely the power ramp-rate, power and energy capability of generators, loads and storage devices, are presented.
The flexibility properties of different power system unit types, e.g. load, generation and storage units that are non-controllable, curtailable or fully controllable are qualitatively analyzed and compared to each other. Quantitative results and flexibility visualizations are presented for intuitive power system examples.
\end{abstract}

\vspace{0.25cm}
%\begin{IEEEkeywords}
%\end{IEEEkeywords}

%\input{Acronym/DEFacronym}
%\chapter*{List of Acronyms}
%\acrodef{}[]{}
%\begin{acronym}

\acrodef{ac}[AC]{Alternating Current}                                 
\acrodef{ami}[AMI]{Automatic Metering Infrastructure}
\acrodef{ams}[AMS]{Asset Management System}
\acrodef{as}[AS]{Ancillary Services}
\acrodef{ase}[ASE]{Aggregated Swing Equation}

\acrodef{bess}[BESS]{Battery Energy Storage System}
\acrodef{bms}[BMS]{Business Management System}

\acrodef{caes}[CAES]{Compressed Air Energy Storage}
\acrodef{ccgt}[CCGT]{Combined Cycle Gas Turbine}
\acrodef{cdd}[CDD]{cooling degree day}
\acrodef{ce}[CE]{Continental European}
\acrodef{cftoc}[CFTOC]{constrained finite-time optimal control}
\acrodef{chp}[CHP]{Combined Heat and Power Plant}
\acrodef{coi}[COI]{Center Of Inertia}
\acrodef{csp}[CSP]{Concentrating Solar Power}

\acrodef{dae}[DAE]{Differential Algebraic Equations}
\acrodef{dc}[DC]{Direct Current}       
\acrodef{dms}[DMS]{Distribution Management System}
\acrodef{dlr}[DLR]{Dynamic Line Rating}
\acrodef{dg}[DG]{Distributed Generation}
\acrodef{dr}[DR]{Demand Response}
\acrodef{dsm}[DSM]{Demand-Side Management}
\acrodef{dso}[DSO]{Distribution System Operator}
\acrodef{dsp}[DSP]{Demand-Side Participation}
\acrodef{dss}[DSS]{Descriptor State Space}

\acrodef{eeg}[EEG]{"Erneuerbare-Energien-Gesetz"}
\acrodef{eex}[EEX]{European Energy Exchange}
\acrodef{eia}[EIA]{Energy Information Agency}
\acrodef{ems}[EMS]{Energy Management System}
\acrodef{entsoe}[ENTSO-E]{European Network of Transmission System Operators for Electricity}
\acrodef{epex}[EPEX]{European Power Exchange}
\acrodef{ev}[EV]{Electric Vehicles}

\acrodef{facts}[FACTS]{Flexible AC Transmission System}
\acrodef{ferc}[FERC]{Federal Energy Regulatory Commission}
\acrodef{fit}[FIT]{Feed-In Tariff}

\acrodef{hdd}[HDD]{heating degree day}
\acrodef{hp}[HP]{Heat Pumps}                                
\acrodef{hpfc}[HPFC]{hourly price forward curve}
\acrodef{hsl}[HSL]{Hydro Storage Lake}
\acrodef{hv}[HV]{High Voltage}                                 
\acrodef{hvac}[HVAC]{Heating, Ventilation, and Air Conditioning}
\acrodef{hvdc}[HVDC]{High Voltage Direct Current}

\acrodef{ice}[ICE]{Internal Combustion Engine}
\acrodef{ict}[ICT]{Information \& Communication Technology}
\acrodef{iea}[IEA]{International Energy Agency}        
\acrodef{ied}[IED]{Intelligent Electronic Device}
\acrodef{ieee}[IEEE]{Institute of Electrical and Electronics Engineers}
\acrodef{iso}[ISO]{Independent System Operator}

\acrodef{lmi}[LMI]{Linear Matrix Inequalities}
\acrodef{lti}[LTI]{Linear Time-Invariant}

\acrodef{mape}[MAPE]{mean absolute prediction error}
\acrodef{mpc}[MPC]{Model Predictive Control}
\acrodef{mse}[MSE]{mean square error}

\acrodef{nis}[NIS]{Network Information System}
\acrodef{nlr}[NLR]{Nominal Line Rating}
\acrodef{npm}[NPM]{Network-Preserving Model}
\acrodef{nreap}[NREAP]{National Renewable Energy Action Plan}
\acrodef{ntc}[NTC]{Net Transfer Capacity}

\acrodef{opf}[OPF]{Optimal Power Flow}
\acrodef{otc}[OTC]{Over-the-Counter}

\acrodef{phev}[PHEV]{Plug-in Hybrid Electric Vehicles}
\acrodef{phs}[PHS]{Pumped Hydro Storage}
\acrodef{pjm}[PJM]{Pennsylvania--New~Jersey--Maryland Interconnection}
\acrodef{pmu}[PMU]{Phasor Measurement Unit}
\acrodef{pshp}[PSHP]{Pumped Storage Hydro Plant}
\acrodef{pst}[PST]{Phase-Shifting Transformer}
\acrodef{purpa}[PURPA]{Public Utility Regulatory Policies Act}
\acrodef{pv}[PV]{Photovoltaics}
\acrodef{pwa}[PWA]{Piece-wise Affine}

\acrodef{ree}[REE]{Red El\'{e}ctrica Espa\~{n}a}
\acrodef{res}[RES]{Renewable Energy Sources}
\acrodef{ror}[ROR]{Run-Of-River Hydro}
\acrodef{rms}[RMS]{Root Mean Square}
\acrodef{rmse}[RMSE]{Root Mean Square Error}
\acrodef{rps}[RPS]{Renewable Portfolio Standards}
\acrodef{rte}[RTE]{R\'{e}seau de Transport d'\'{E}lectricit\'{e}}

\acrodef{soc}[SOC]{State-of-Charge}
\acrodef{scada}[SCADA]{Supervisory Control And Data Acquisition}

\acrodef{td}[T\&D]{Transmission and Distribution}
\acrodef{tso}[TSO]{Transmission System Operator}

\acrodef{vpp}[VPP]{Virtual Power Plants}

\acrodef{wampac}[WAMPAC]{Wide-Area Monitoring, Protection And Control}

%\end{acronym}

%\IEEEpeerreviewmaketitle

\maketitle

\section{Introduction}

This paper presents a novel approach for analyzing the available operational flexibility of a given power system. In the context of this paper we mean by this the combined available operational flexibility that an ensemble of diverse power system units in a geographically confined grid zone can provide in each time-step during the operational planning, given load demand and \acf{res} forecast information, as well as in real-time in case of a contingency. Operational flexibility is essential for mitigating disturbances in a power system such as outages or forecast deviations of either power in-feed, i.e. from wind turbines or solar units, or power out-feed, i.e. load demand. Metrics for assessing the technical operational flexibility of power systems, i.e. power ramp-rate ($\rho$), power capacity ($\pi$) and energy capacity ($\epsilon$), have been proposed by Makarov et al. in~\cite{Makarov:2009} and their meaning further discussed by the authors in~\cite{Ulbig:2012}. 
In this paper we establish the necessary framework for quantifying and visualizing the technically available operational flexibility of individual power system units and ensembles thereof. The functional modeling of all power system units is accomplished using the Power Nodes modeling framework introduced in~\cite{PowerNodes:2010, PowerNodes:2012}.
The flexibility properties of different power system unit types, e.g. load, generation and storage units that are non-controllable, curtailable or fully controllable are qualitatively analyzed and compared to each other.
Quantitative results as well as flexibility visualizations of the here proposed flexibility assessment framework are presented for intuitive benchmark power systems.

The remainder of this paper is organized as follows: Section~\ref{sec:OpFlex} discusses operational flexibility and its role in power system operation. It also introduces necessary metrics for operational flexibility. Section~\ref{sec:OpFlexModel} explains how operational flexibility can be modeled using the Power Nodes functional modeling framework. This is followed by Section~\ref{sec:OpFlexAnalysis}, which illustrates how operational flexibility can be quantified and analyzed for individual power system units as well as for unit ensembles. Finally, a conclusion and a summary of the contributions of this paper are given in Section~\ref{sec:Conclusion}.

\section{Operational Flexibility in Power Systems}\label{sec:OpFlex}

Operational flexibility is an important property of electric power systems and essential for mitigating disturbances in a power system such as outages or forecast deviations of either power in-feed, i.e. from wind turbines or \acf{pv} units, or power out-feed, i.e. load demand. The availability of sufficient operational flexibility is a necessary prerequisite for the effective grid integration of large shares of fluctuating power in-feed from variable \ac{res}.

\subsection{Increasing Need for Operational Flexibility}
In recent years power system dispatch optimization and real-time operation are more and more driven by several major trends which include notably 
\begin{enumerate}
	\item Wide-spread deployment of variable \ac{res}, i.e.~wind turbines and \ac{pv} units, has led to significant relative and absolute shares of power generation that is highly fluctuating and not perfectly predictable nor fully controllable. Variable \ac{res} power in-feed causes non-deterministic power imbalances and power flow changes on all grid levels~\cite{IEA2005, Jones:2014}.
	\item Growing power market activity has led to operational concerns of its own, i.e.~deterministic frequency deviations caused by transient power imbalances due to more frequent changes in the now market-driven operating set-point schedules of power plants as well as more volatile (cross-border) power flow patterns~\cite{Weissbach:2008}.

%\afterpage{\restoregeometry}
%\restoregeometry

\begin{figure}[t]
  %\vspace{-0.2cm}
  \centering               %left bottom right top
  {\includegraphics[trim = 0cm   0cm   0cm   0cm, clip=false, angle=0, width=0.80\linewidth, keepaspectratio, draft=false]{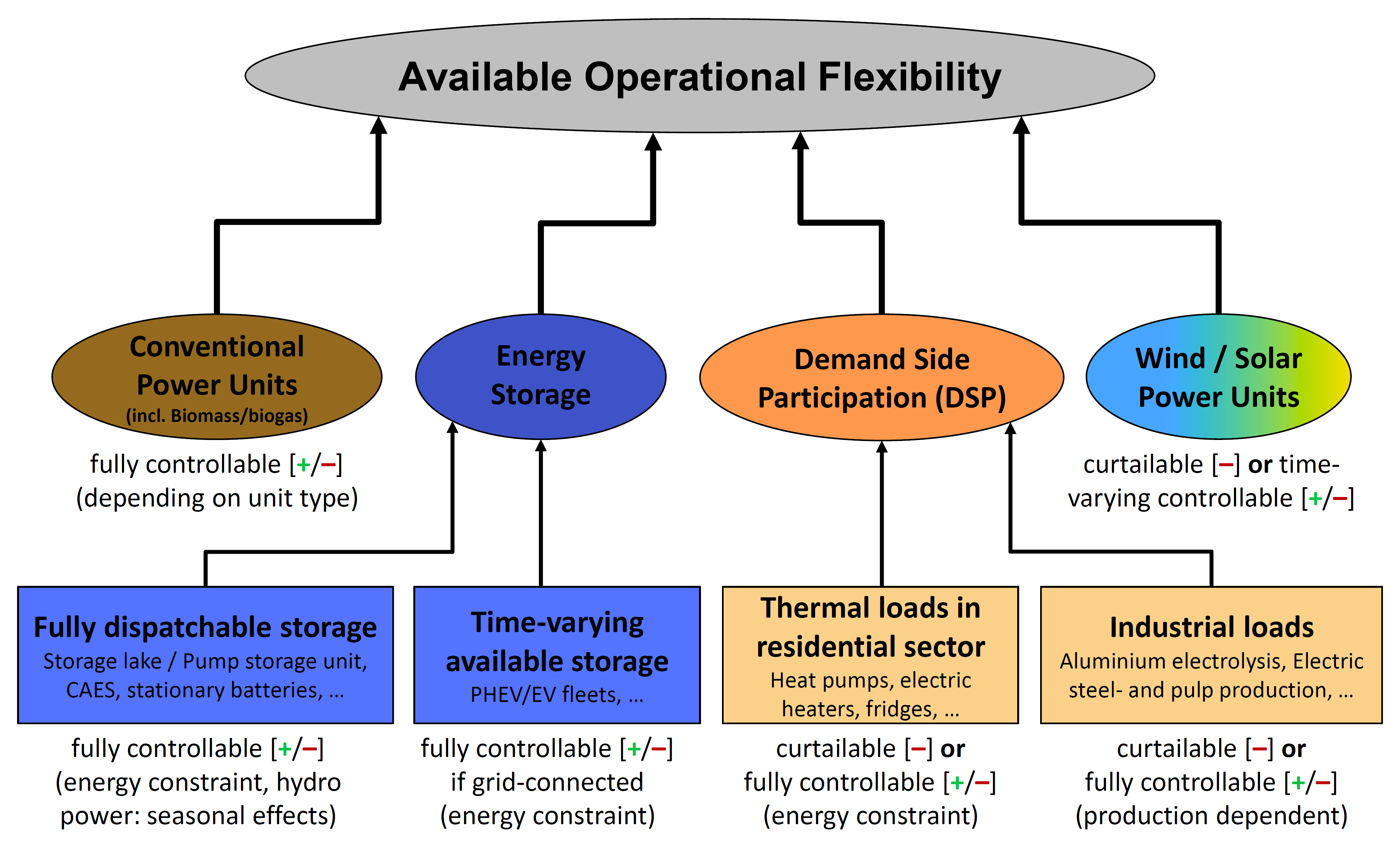}}
  \vspace{-0.15cm}
  \caption{Sources of Operational Flexibility in Power Systems.} 
  \label{fig:FlexSources}
  \vspace{-0.250cm}
\end{figure}
  
	\item The emergence of a \emph{smart grid} notion or vision as a driver for change in power system operation~\cite{Smartgrid}. Using the reference framework of control theory, the term \emph{smart grid} can be understood as the sum of all efforts that improve observability and controllability over individual power system processes, i.e.~power in-feed to the grid and power
out-feed from the grid as well as power flows on the demand/supply side, happening on all voltage levels of the electricity
grid. An improved observability and controllability of individual power system units should also lead to an improved
observability and controllability of the entire power system and the processes happening therein.	
\end{enumerate}
 
Altogether, these developments constitute a major paradigm shift for the management of power systems. Operating power systems optimally in this more complex environment requires a more detailed assessment of \emph{available} operational flexibility at every point in time for effectively mitigating the outlined disturbances. 
 
Operation flexibility in power system operation and dispatch planning is of importance and has a significant commercial value. Ancillary service markets enable system operators the cost-effective procurement of needed control reserve products. In the case of frequency control schemes which are in essence a set of differently structured flexibility services provided to system operators for achieving active power regulation on different time-scales~\cite{Kundur1994}, the overall remuneration for providing control power and energy on ancillary service markets is usually significantly higher than for bulk energy from spot markets~\cite{Regelleistung2013}. The value of operational flexibility can also be shown indirectly by looking at the inflexibility costs incurred by conventional generation units in the form of ramping costs as well as power plant start/stop costs. 
In some power markets, the real or merely perceived inflexibility of generator units to reduce their power output from planned set-points appears in the form of negative bids in the supply-side curve of the merit order~\cite{EPEX, Nicolosi:2010}. Negative bids may either reflect costs that would be incurred in case a plant's power output is lowered, e.g. lower efficiency as well as wear and tear, or the goal to keep a certain power plant online, i.e.~must-run units that provide ancillary services or \ac{res} units that have in-feed priority.

\subsection{Sources of Operational Flexibility}
Different sources of power system flexibility exist as is illustrated in Fig.~\ref{fig:FlexSources}. 
Operational flexibility can be obtained on the generation-side in the form of dynamically fast responding conventional power plants, e.g.~gas or oil-fueled turbines or rather flexible modern coal-fired power plants and on the demand-side by means of adapting the load demand curve to partially absorb fluctuating RES power in-feed. In addition to this, \ac{res} power in-feed can also be curtailed or, in more general terms, modulated below its given time-variant maximum output level. Furthermore, stationary storage capacities, e.g. hydro storage, \ac{caes}, stationary battery or fly-wheel systems, as well as time-variant storage capacities, e.g.~electric vehicle fleets, are well-suited for providing operational flexibility.

Additional flexibility can be obtained from other grid zones via the electricity grid's tie-lines in case that the available operational flexibility in one's own grid zone is not sufficient or more expensive than elsewhere. Power import and export, nowadays facilitated by more and more integrated transnational power markets, is used in daily power system operation to a certain degree as a \emph{slack bus} for fulfilling the active power balance and mitigating power flow problems of individual grid zones by tapping into the flexibility potential of other grid zones. For power system operation, importing needed power in certain situations and exporting undesirable power in-feed in other situations to neighboring grid zones is for the time being probably the most convenient and cheapest measure for increasing operational flexibility. However, power import/export can only be performed within the limits given by the agreed line transfer capacities between the grid zones. In the European context this corresponds to the \ac{ntc} values~\cite{ENTSOE:2011}, which are a rather conservative measure of available grid electricity transfer capacity.

In liberalized power systems, operational flexibility is traded in the form of \emph{energy products} via power markets, i.e.~day-ahead and intra-day spot markets, as well as \emph{control reserve products}, i.e.~primary/secondary/tertiary frequency control reserves, from \ac{as} markets.
 
\subsection{Definitions of Operational Flexibility}

%\begin{itemize}
%	\item Flexibility Definition~\cite{Makarov:2009}
%	\item Power Nodes~\cite{PowerNodes:2010, PowerNodes:2012}
%	\item Operational Flexibility~\cite{Ulbig2012}
%\end{itemize}

The term Operational Flexibility in power systems, or simply flexibility is often not properly defined and may refer to very different things, ranging from the quick response times of certain generation units, e.g. gas turbines, to the degree of efficiency and robustness of a given power market setup. The topic has received wide attention in recent years~\cite{Makarov:2009, Kirschen:2012, OMalley:2012, Ulbig:2012, Huber:2014, Jones:2014}.

In the following the focus is on the basic \emph{technical} capability of individual power system units to modulate power and energy in-feed into the grid, respectively power out-feed out of the grid.

\subsection{Metrics for Operational Flexibility}

\begin{figure}[bht]
\vspace{-0.5cm}
%\begin{picture}(80,140)
%\put(0,0)
\centering               %left bottom right top
{\includegraphics[trim = 0cm   0cm   0cm   0cm, clip=true, angle=0, width=0.75\linewidth, keepaspectratio, draft=false]{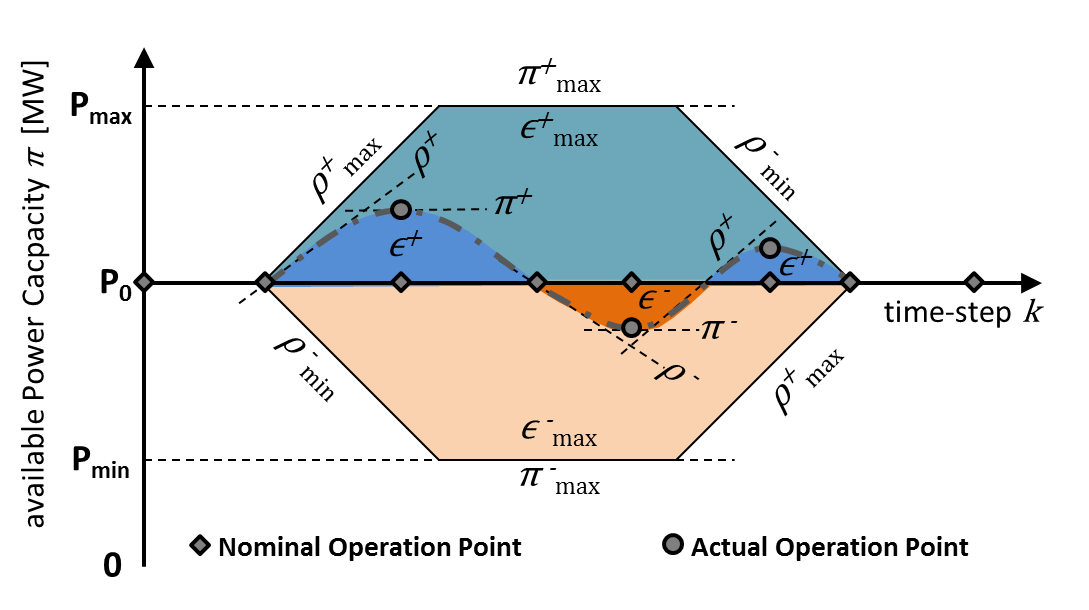}}
%\end{picture}
\caption{Flexibility Metrics in Power Systems Operation: \newline Power Ramp-Rate $\rho$, Power $\pi$ and Energy $\epsilon$.} \label{fig:Flex_Defs}
\vspace{-0.25cm}
\end{figure}

For analysis purposes, this technical capability needs to be characterized and categorized by appropriate flexibility metrics. A valuable method for assessing the needed operational flexibility of power systems, for example for accommodating high shares of wind power in-feed, has been proposed by Makarov et al. in~\cite{Makarov:2009}. There, the following metrics have been characterized:
\begin{itemize}
	\item Power provision capacity $\pi$ (MW),
	\item Power ramp-rate capacity $\rho$ (MW/min.),
	\item Energy provision capacity $\epsilon$ (MWh) as well as 
	\item Ramp duration $\delta$ (min.).
\end{itemize}
%\vspace{+0.2cm}
Their role in modulating the operation point of a power system unit and with it the relative power flow into the grid ($>0$) and out of the grid ($<0$) with respect to the original operation point is depicted in Fig.~\ref{fig:Flex_Defs}. % Note that the term $\delta$ is actually dependent on $\pi$ and $\rho$ (and vice versa); $\delta = \frac{\pi}{\rho}$.
Here, the deliberate deviation between the nominal power plant output trajectory and the actual power output trajectory is bounded by the maximum flexibility capability, i.e.~the three metrics $\rho$, $\pi$ and $\epsilon$, of the power plant in question. For the sake of simplicity and clarity we will stick to the same notation as in~\cite{Makarov:2009}.

Having a closer look on the proposed flexibility metrics, the following two things can be observed:
\begin{itemize}
	\item The ramp duration $\delta$ is actually dependent on the power ramp rate $\rho$ and power capacity $\pi$ (and vice versa) as $\delta = \pi/\rho$.
	\item It it thus entirely sufficient to focus on the power-related metrics $\rho$, $\pi$ and $\epsilon$ for describing operational flexibility.
	\item An intriguing feature is that the metric terms $\rho$, $\pi$ and $\epsilon$ are closely linked via integration and differentiation operations in the time domain. The interaction of the individual metrics clearly exhibit so-called double integrator dynamics: energy is the integral of power, which in turn is the integral of power ramp-rate~(Eq.~\ref{equ:Trinity}). The three metrics constitute a \emph{flexibility trinity} in power system operation, as they cannot be thought of independently due to the inter-temporal linking. 
\end{itemize}
\begin{IEEEeqnarray}{rCL}\label{equ:Trinity}
\left.
\setlength{\arraycolsep}{3pt}
\begin{array}{lclclclclcl}
%&\textbf{Ramp-Rate}&  &\frac{d}{dt}&      &\textbf{Power}&  &\frac{d}{dt}&      &\textbf{Energy}&   \\
%&&                    &&                  &&                &&                  &&                  \\
%&\mathbf{\rho}&       &\leftrightarrows&  &\mathbf{\pi}&    &\leftrightarrows&  &\mathbf{\epsilon}& \\
%&&                    &&                  &&                &&                  &&                  \\
%&\textrm{[MW/min.]}&    &\int dt&           &\textrm{[MW]}&     &\int dt&           &\textrm{[MWh]}&      \\
%&&                    &\dfrac{d}{dt}&     &&                &\dfrac{d}{dt}&     &&                  \\
%&\mathbf{\rho}&       &\leftrightarrows&  &\mathbf{\pi}&    &\leftrightarrows&  &\mathbf{\epsilon}& \\
%&&                    &\int dt&           &&                &\int dt&           &&                  \\
%&&                    &&                  &&                &&                  &&                  \\
%&\textbf{Ramp-Rate}&  &&                  &\textbf{Power}&  &&                  &\textbf{Energy}&   \\
%&\textrm{[MW/min.]}&    &&                  &\textrm{[MW]}&     &&                  &\textrm{[MWh]}&      \\
&\textbf{Ramp-Rate}&  &&                  &\textbf{Power}&  &&                  &\textbf{Energy}&   \\
&\textrm{[MW/min.]}&    &&                  &\textrm{[MW]}&     &&                  &\textrm{[MWh]}&      \\
%&&                    &&                  &&                &&                  &&                  \\
%&&                    &\int dt&           &&                &\int dt&           &&                  \\
%&&                    &\Rightarrow&       &&                &\Rightarrow&       &&                  \\
%&\mathbf{\rho}&       &&                  &\mathbf{\pi}&    &&                  &\mathbf{\epsilon}& \\
%&&                    &&                  &&                &&                  &&                  \\
%&&                    &\Leftarrow&        &&                &\Leftarrow&        &&                  \\
%&&                    &\frac{d}{dt}&      &&                &\frac{d}{dt}&      &&                  \\
&&                    &\int dt&           &&                &\int dt&           &&                  
\\
&\mathbf{\rho}&       &\rightleftarrows&  &\mathbf{\pi}&    &\rightleftarrows&  &\mathbf{\epsilon}& \\
&&                    &\frac{d}{dt}&      &&                &\frac{d}{dt}&     &&                  
\\
&&                    &&                  &&                &&                  &&                  \\
&&                    &&                  &&                &&                  &&                  \\
&&                    &&                  &&                &&                  &&                  \\
\end{array}
\right.
\vspace{-1.00cm}
\end{IEEEeqnarray}
Using these three flexibility metrics instead of only one, for example the power ramping capability $\rho$ as in~\cite{OMalley:2012}, allows a more accurate and complete representation of power system flexibility, including the relevant inter-temporal constraints over a given time interval. The power ramp-rate for absorbing a disturbance event, measured in MW/min, in a power system may be abundant at a certain time instant. But for a persistent disturbance, the maximum regulation power that can be provided by a generator is limited as is the maximum regulation energy that can be provided by storage units, which are inherently energy-constrained. As the share of storage units in power systems and their importance for the grid integration of \ac{res} in-feed is rising, the inter-temporal links between providing ramping capability and eventually reaching power/energy limits cannot be neglected when assessing the available operational flexibility of a power system. 

\begin{figure}[htb!]
%\vspace{-1.0cm}
%\begin{picture}(80,140)
%\put(0,0)
\centering               %left bottom right top
{\includegraphics[trim = 0cm   0cm   0cm   0cm, clip=true, angle=0, width=0.5\linewidth, keepaspectratio, draft=false]{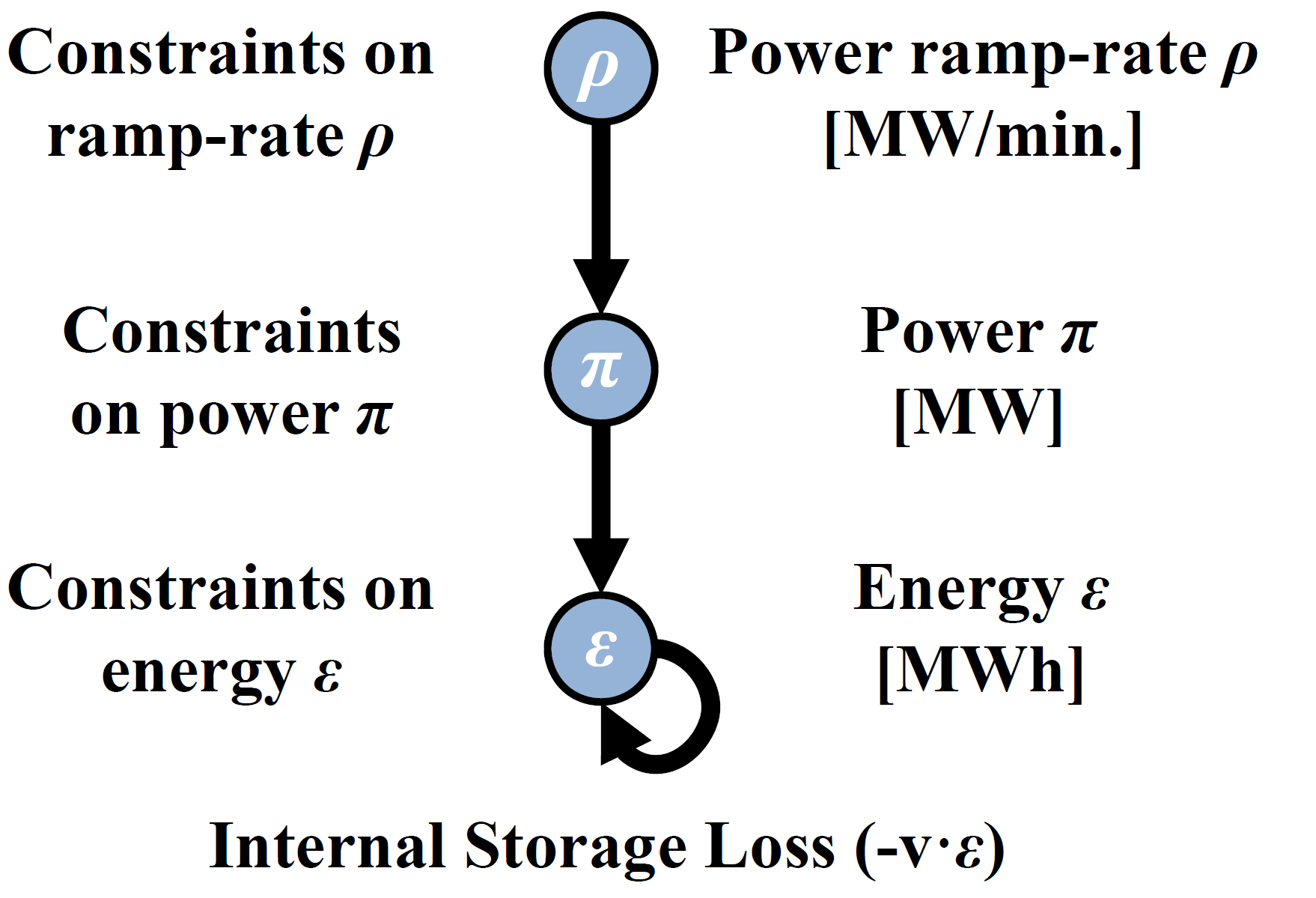}}
%\end{picture}
\caption{Inter-temporal linking of flexibility metrics including internal storage losses (dissipation).} \label{fig:Flex_Link}
\vspace{-0.25cm}
\end{figure}

Having defined these flexibility metrics as well as the causal inter-linking between them~(Fig.~\ref{fig:Flex_Link}), allows the assessment of the available operational flexibility of an individual power system unit and for whole power systems. Note that the operational constraints, i.e.~min/max ramping, power and energy constraints, of individual power system units have to be considered when assessing their available operational flexibility.

\section{Modeling of Operational Flexibility}\label{sec:OpFlexModel}

\begin{figure}[htb!]
%\vspace{-1.0cm}
%\begin{picture}(80,140)
%\put(0,0)
\centering               %left bottom right top
{\includegraphics[trim = 0cm   0cm   0cm   0cm, clip=true, angle=0, width=0.60\linewidth, keepaspectratio, draft=false]{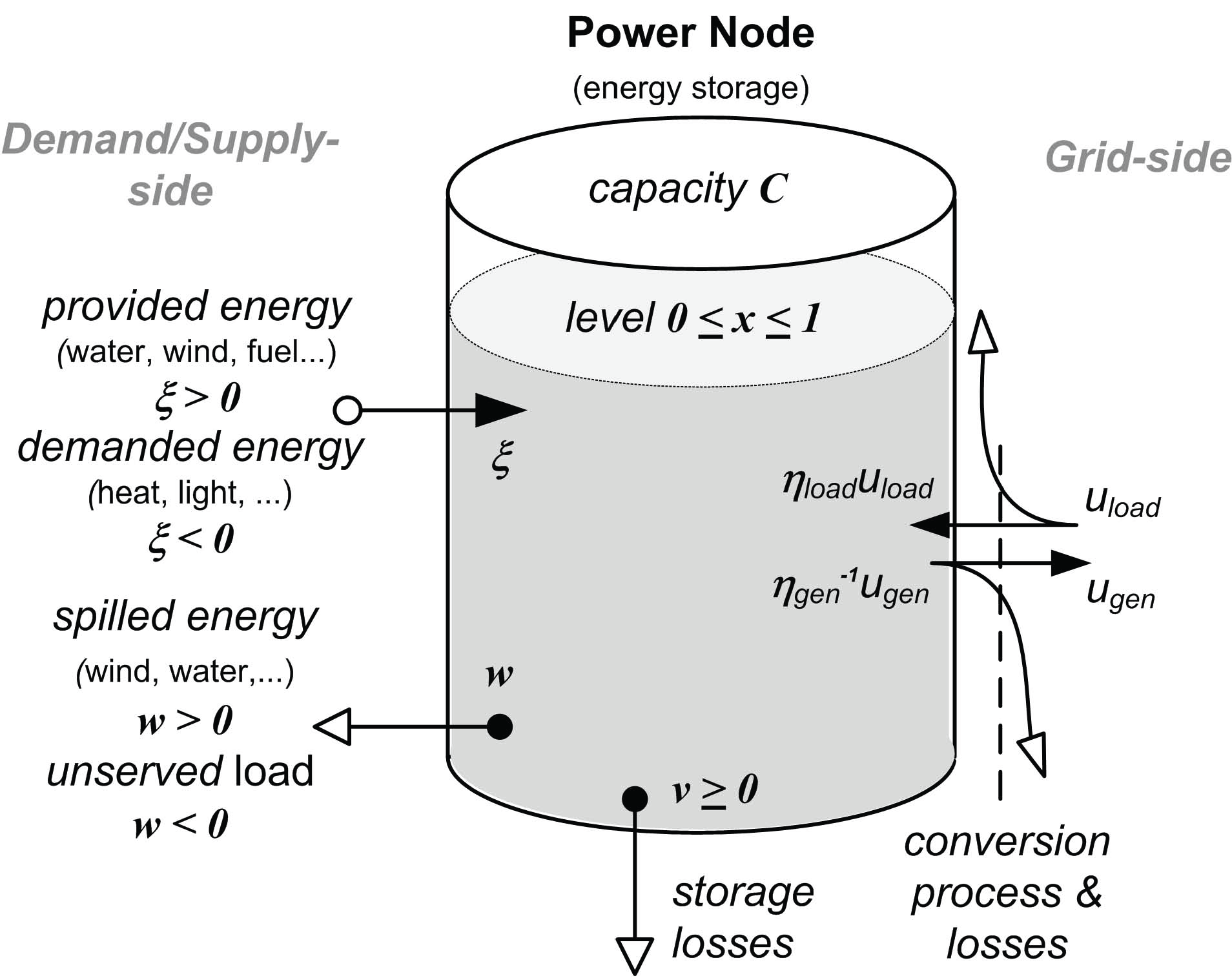}}
%\end{picture}
\caption{Power Node model of an energy storage unit with power in-feed ($u_{\textrm{gen}}$) and out-feed ($u_{\textrm{load}}$).} \label{fig:PN}
\vspace{-0.25cm}
\end{figure}

The analysis and assessment of operational flexibility first of all necessitates a modeling framework that allows to explicitly include information on the degree of freedom for shifting operation set-points so as to modulate the power in-feed and out-feed patterns of individual power system units. This includes information on whether or not a unit has a storage and is thus energy-constrained, whether or not a unit provides fluctuating power in-feed, and what type of controllability and observability, including predictability that a system operator has over fluctuating generation and demand processes (i.e.~full, partial or none). The combination of all these properties defines a unit's operational flexibility.

For our modeling purposes we use the Power Nodes modeling framework, which allows the detailed functional modeling of power system units such as
%\textbf{XXX Include Paper Citations XXX}
\begin{itemize}
  \item diverse storage units, e.g.~batteries, fly-wheels, pumped hydro, \ac{caes}, ...,
  \item diverse generation units, e.g.~fully dispatchable conventional generators, variably in-feeding power units, e.g.~wind turbine and PV units, and
  \item diverse load units, e.g.~conventional (non-controllable), interruptible or thermal (both partially controllable), ...,
\end{itemize}
including their operational constraints as well as relevant information of their underlying power supply and demand processes. Operational constraints such as min/max ramp rates, min/max power set-points and energy storage operation ranges, information of the underlying power system processes (i.e.~fully controllable, curtailable/sheddable or non-controllable) as well as information on observability and predictability of underlying power system processes (i.e.~state measures and/or state-estimation and prediction of fully or only partially observable/predictable system and control input states) can also be included. 
The workings of the Power Node notation are illustrated by the model representation of an energy storage unit~(Fig.~\ref{fig:PN}). The provided and demanded energies are lumped into an external process termed $\xi$, with $\xi < 0$ denoting energy use and $\xi > 0$ energy supply. The term $u_{\textrm{gen}}$ describes a conversion corresponding to a power generation with an efficiency $\eta_{\textrm{gen}}$, while $u_{\textrm{load}}$ describes a conversion corresponding to consumption with an efficiency $\eta_{\textrm{load}}$. The introduction of generic energy storages in the Power Nodes framework adds a modeling layer to classical power system modeling. Its energy storage level, the State-of-Charge (SOC), is normalized to $0 \le x \le 1$ with an energy storage capacity $C \ge 0$. The illustrated storage unit serves as a buffer between the external process $\xi$ and the two grid-related power exchanges $u_{\textrm{gen}} \ge 0$ and $u_{\textrm{load}} \ge 0$. Internal energy losses associated with energy storage, e.g.~physical, state-dependent dissipation losses, are modeled by the power dissipation term $v(x) \ge 0$, while enforced energy losses, e.g.~curtailment/shedding of a power supply or demand process, are denoted by the waste power term $w$, where $w > 0$ denotes a loss of provided energy and $w < 0$ an unserved load demand. 

The dynamics of a power node $i \in \mathcal{N} = \{1, \dots, N\}$, which can be nonlinear in the general case, are:
\begin{eqnarray}\label{eq:general_powernode}
	C_i \, \dot{x}_i &=& \eta_{\textrm{load},i}  \, u_{\textrm{load},i} - \eta^{-1}_{\textrm{gen},i}  \, u_{\textrm{gen},i} +\, \xi_i - w_i - v_i, \nonumber \\
\textrm{s.t.} &\textrm{(a)}&  0 \leq {x}_i \leq 1 \quad , \nonumber \\
	&\textrm{(b)}&  0 \leq u_{\textrm{gen},i}^{\textrm{min}} \leq u_{\textrm{gen},i} \leq u_{\textrm{gen},i}^{\textrm{max}} \nonumber \quad,\\
	&\textrm{(c)}&  0 \leq u_{\textrm{load},i}^{\textrm{min}} \leq u_{\textrm{load},i} \leq u_{\textrm{load},i}^{\textrm{max}} \nonumber \quad,\\
	&\textrm{(d)}& \dot u_{\textrm{gen},i}^{\textrm{min}} \leq \dot u_{\textrm{gen},i} \leq \dot u_{\textrm{gen},i}^{\textrm{max}} \nonumber \quad, \\
	&\textrm{(e)}& \dot u_{\textrm{load},i}^{\textrm{min}} \leq \dot u_{\textrm{load},i} \leq \dot u_{\textrm{load},i}^{\textrm{max}} \nonumber \quad, \\
    &\textrm{(f)}& 0 \le \xi_i \cdot w_i  \quad, \nonumber \\
    &\textrm{(g)}& 0 \le |\xi_i| - |w_i| \nonumber \quad, \\
    &\textrm{(h)}& 0 \le v_i    \quad.
%    &\textrm{(h)}& 0 \le v_i    \quad \forall \, i=1, \dots, N  \nonumber \quad.
\end{eqnarray}
Depending on the specific process represented by a Power Node, each term in the Power Node equation may be controllable or not, observable or not, and driven by an external process or not. Internal dependencies, such as a state-dependent loss term vi(xi), are possible. Charge and discharge efficiencies may be non-constant and possibly also state-dependent: $\eta_{\textrm{load,i}}(x_i)$, $\eta_{\textrm{gen,i}}(x_i)$. Non-linear conversion efficiencies can be arbitrarily well approximated by a set of \ac{pwa} linear equations. The constraints (a)--(h) denote a generic set of requirements on the variables. They are to express that (a) the state of charge is normalized, (b)--(e) the grid power in-feeds and out-feeds as well as their time derivatives (ramp-rates) are non-negative and constrained, (f) the supply or demand and the curtailment need to have the same sign, (g) the supply/demand curtailment cannot exceed the supply/demand itself, and (h) the storage losses are non-negative.
The explicit mathematical form of a power node equation depends on the particular modeling case. The notation provides technology-independent categories that can be linked to evaluation functions for energy and power balances. Power nodes can also represent energy processes that are independent of storage, such as fluctuating \ac{res} generation. 
%A process without storage implies an algebraic coupling between the instantaneous quantities $\xi_i$, $w_i$, $u_{\textrm{gen},i}$, and $u_{\textrm{load},i}$. Storage-dependent losses do not exist in this case ($v_i=0$). Equation~(\ref{eq:general_powernode}) thus degenerates to the algebraic constraint
%\begin{equation}
%\xi_i - w_i  = \eta^{-1}_{\textrm{gen},i}\, u_{\textrm{gen},i} - \eta_{\textrm{load},i} \, u_{\textrm{load},i} \quad .%,
%\label{eq:no_storage}
%\end{equation}

More details on the Power Node modeling framework, modeling examples and reasoning can be obtained from~\cite{PowerNodes:2010, PowerNodes:2012}.

%\textbf{XXXXXXXXXXXXXXXXXXXXXXXXX}
\section{Analyzing Operational Flexibility}\label{sec:OpFlexAnalysis}

The functional representation of complex power system interactions using the Power Nodes notation allows a straight-forward analysis of the three power-related operational flexibility metrics, i.e.~power ramp-rate $\rho$, power $\pi$ and energy capability $\epsilon$.

\subsection{Quantification of Operational Flexibility}

\begin{figure}[htb!]
%\vspace{-1.0cm}
%\begin{picture}(80,140)
%\put(0,0)
\centering               %left bottom right top
{\includegraphics[trim = 0cm   0cm   0cm   0cm, clip=true, angle=0, width=0.75\linewidth, keepaspectratio, draft=false]{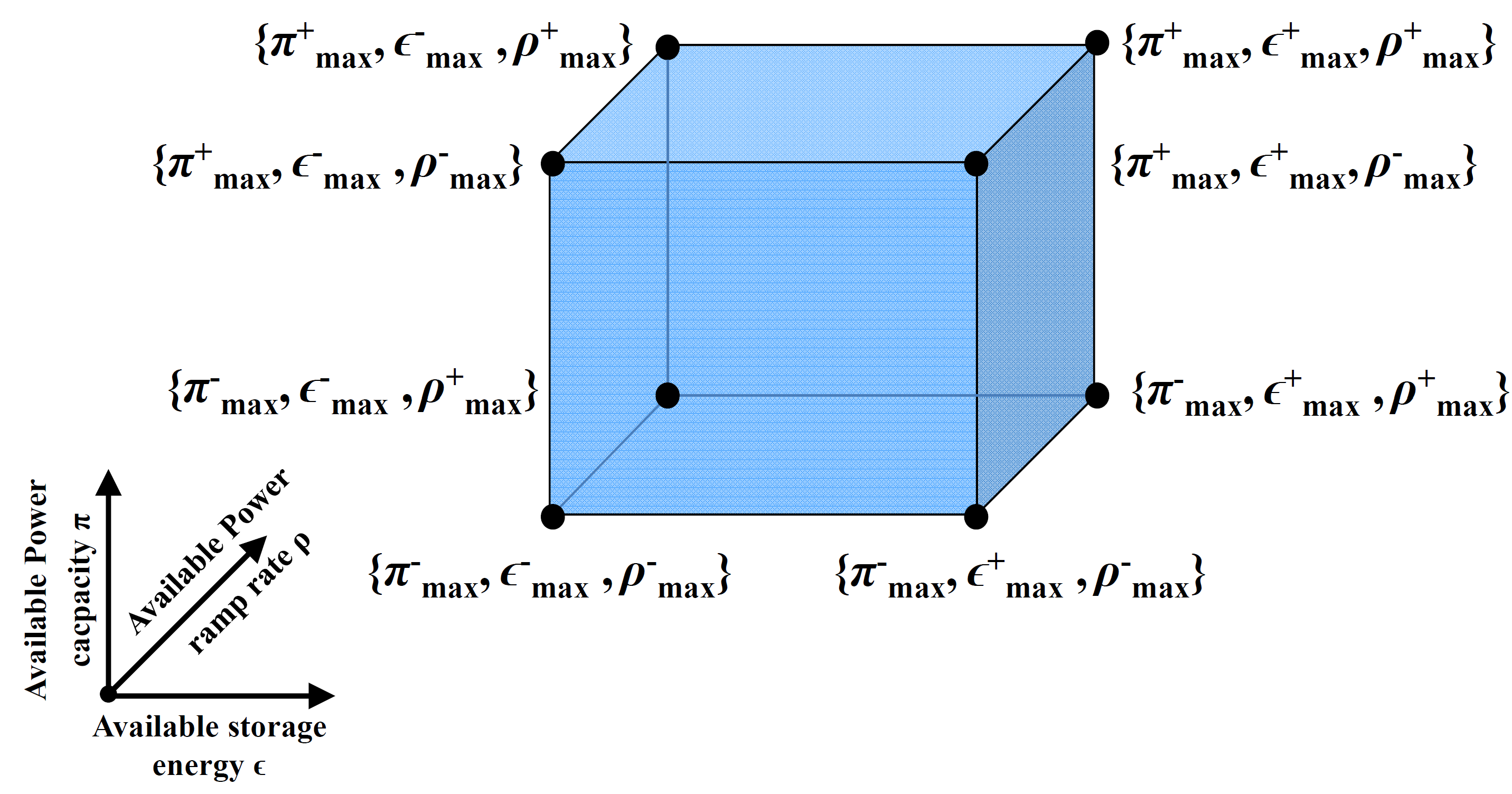}}
%\end{picture}
\caption{Flexibility cube of maximum available operational flexibility of a generic power system unit.} \label{fig:Flex_Cube}
\vspace{-0.25cm}
\end{figure}

\begin{figure}[htb!]
%\vspace{-1.0cm}
%\begin{picture}(80,140)
%\put(0,0)
\centering               %left bottom right top
%{\includegraphics[trim = 1.5cm   1cm   1.5cm   1cm, clip=true, angle=0, width=0.325\linewidth, keepaspectratio, draft=false]{Figure6_Evolution_k1_HQ.png}}
%{\includegraphics[trim = 1.5cm   1cm   1.5cm   1cm, clip=true, angle=0, width=0.325\linewidth, keepaspectratio, draft=false]{Figure6_Evolution_k12_HQ.png}}
%{\includegraphics[trim = 1.5cm   1cm   1.5cm   1cm, clip=true, angle=0, width=0.325\linewidth, keepaspectratio, draft=false]{Figure6_Evolution_k24_HQ.png}}
%{\includegraphics[trim = 1.5cm   1cm   1.5cm   1cm, clip=true, angle=0, width=0.325\linewidth, keepaspectratio, draft=false]{Figure6_Evolution_k36_HQ.png}}
%{\includegraphics[trim = 1.5cm   1cm   1.5cm   1cm, clip=true, angle=0, width=0.325\linewidth, keepaspectratio, draft=false]{Figure6_Evolution_k48_HQ.png}}
%{\includegraphics[trim = 1.5cm   1cm   1.5cm   1cm, clip=true, angle=0, width=0.325\linewidth, keepaspectratio, draft=false]{Figure6_Evolution_k60_HQ.png}}
%\caption{Time-evolution of maximum available operational flexibility (k = 1h, 12h, 24h, 36h, 48h, 60h).} \label{fig:Flex_Evolution}
{\includegraphics[trim = 0.5cm  0.5cm  0.5cm  2cm, clip=true, angle=0, width=0.4\linewidth, keepaspectratio, draft=false]{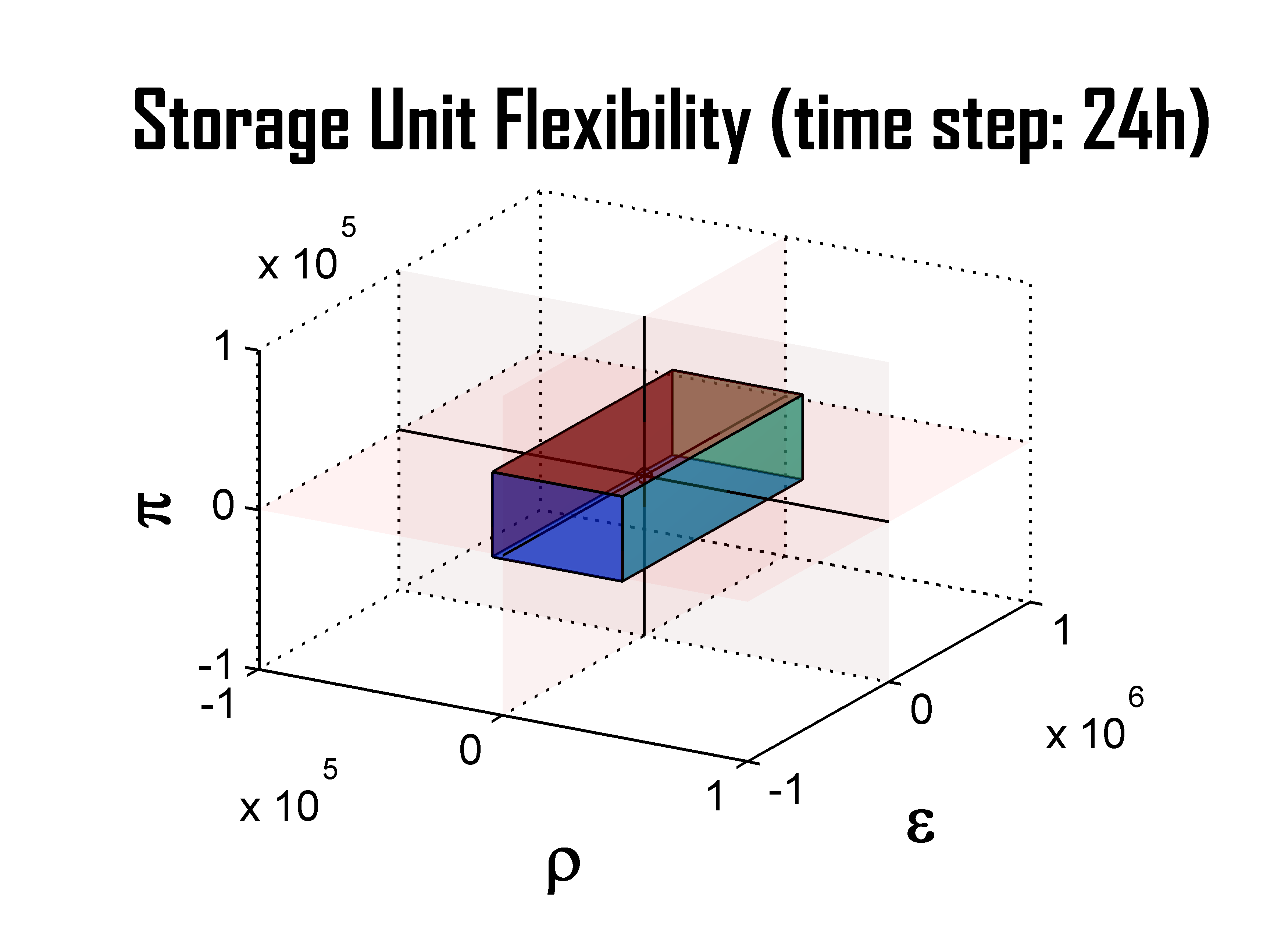}}
{\includegraphics[trim = 0.5cm  0.5cm  0.5cm  2cm, clip=true, angle=0, width=0.4\linewidth, keepaspectratio, draft=false]{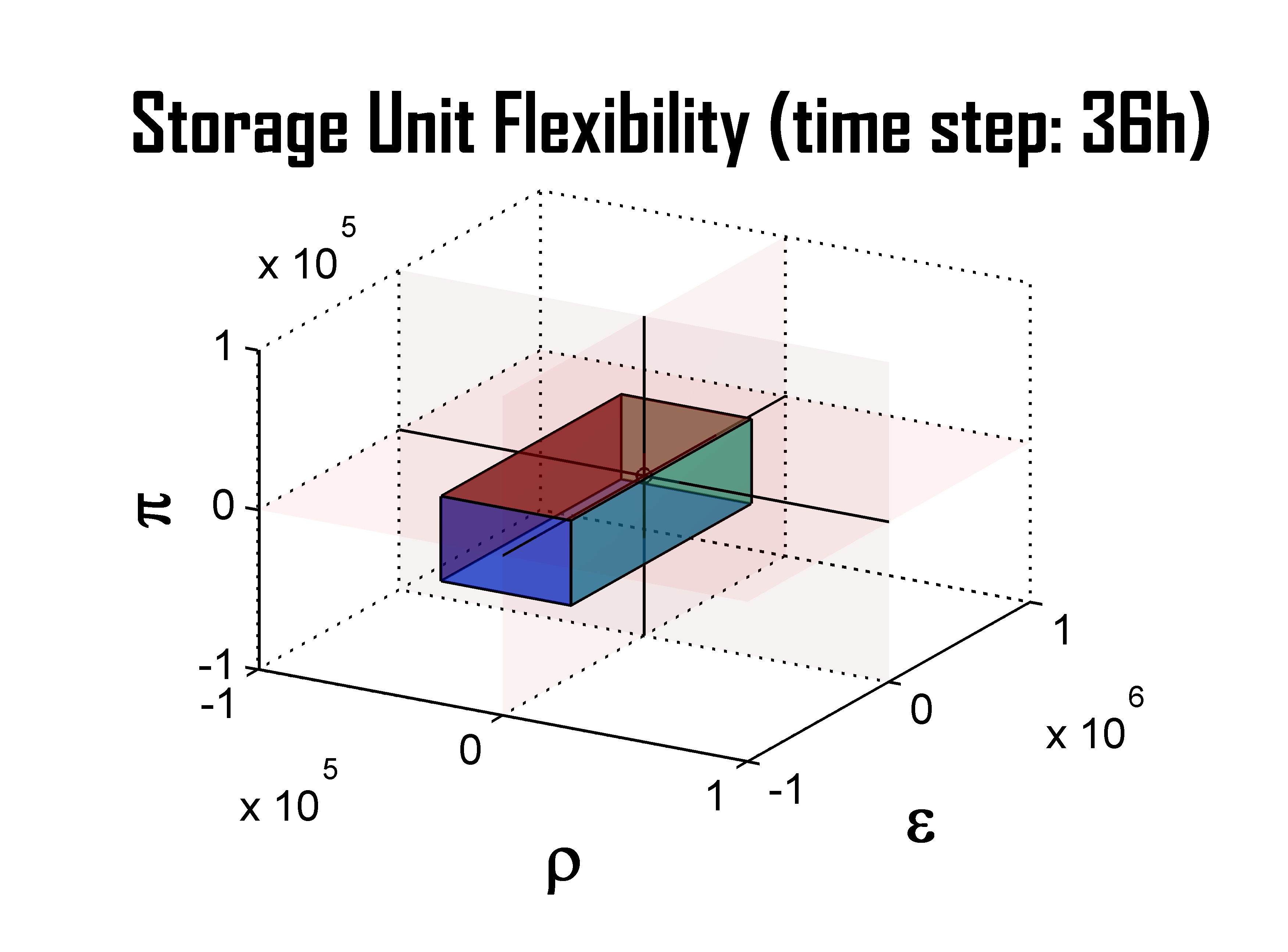}}
{\includegraphics[trim = 0.5cm  0.5cm  0.5cm  2cm, clip=true, angle=0, width=0.4\linewidth, keepaspectratio, draft=false]{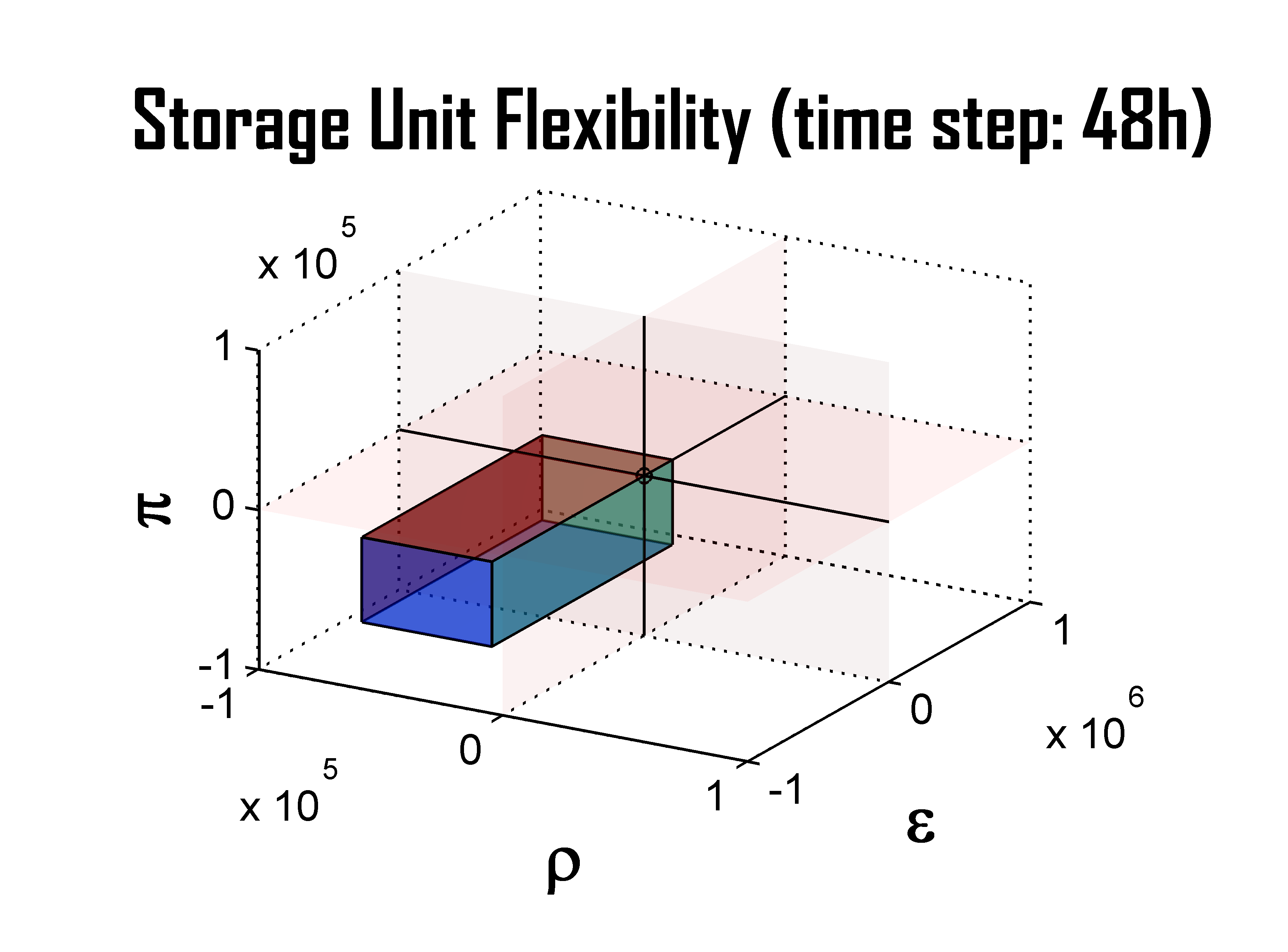}}
{\includegraphics[trim = 0.5cm  0.5cm  0.5cm  2cm, clip=true, angle=0, width=0.4\linewidth, keepaspectratio, draft=false]{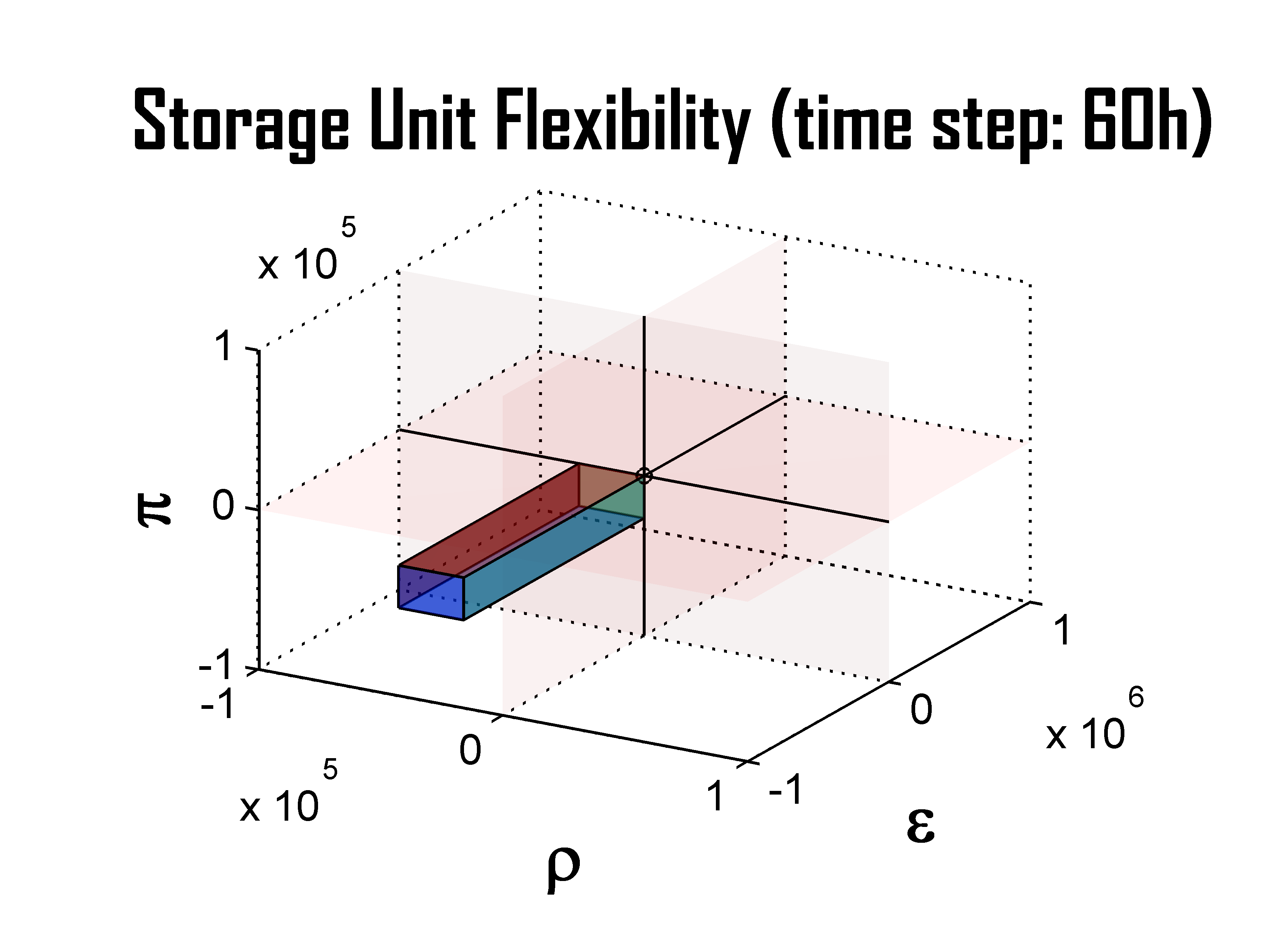}}
%\end{picture}
\caption{Time-evolution of maximum available operational flexibility ($k = 24\,\text{h},\; 36\,\text{h},\; 48\,\text{h},\; 60\,\text{h}$).} \label{fig:Flex_Evolution}
\vspace{-0.25cm}
\end{figure}

Taking as an example the operational flexibility of a generation unit $i$ that has an inherent storage function and the possibility for curtailment, e.g. a hydro storage lake, given by the Power Node model
\begin{equation}
C_i \, \dot{x}_i = - \eta^{-1}_{\textrm{gen},i}  \, u_{\textrm{gen},i} +\, \xi_i - w_i-v_i \quad ,
\end{equation}
for providing power regulation is accomplished by calculating the set of all feasible power regulation points~$\{\pi_i^{\pm}(k)\}$, where up/down power regulation is denoted by~'$+/-$' respectively, based on equation
\begin{IEEEeqnarray}{ccl}\label{eq:FlexAssessmentPi1}
    \Big\{\pi_i^{\pm}(k)\Big\} &=& \Big\{u_{\textrm{gen},i}^{\textrm{feasible}}(k)\Big\} - u_{\textrm{gen},i}^0(k)  \ \ \ \\
                 %&=& u_{\textrm{gen}}^{\frac{\textrm{max.}}{\textrm{min.}}}(k) \ldots \nonumber \ \ \ \\
                 &=& \Big\{\eta_{\textrm{gen}} \cdot \left( \xi - w - v_x - C \dot{x} \right) \Big\}_{k,i} - u_{\textrm{gen},i}^0(k) \nonumber \\
    \textrm{s.t.} && \ 0 \leq u_{\textrm{gen},i}^{\textrm{min}}(k) \leq \{u_{\textrm{gen}}^{\textrm{feasible}}(k)\} \leq u_{\textrm{gen},i}^{\textrm{max}}(k) \nonumber \quad .
    %w^{\frac{\textrm{min}}{\textrm{max}}}
\end{IEEEeqnarray}
Here, $u_{\textrm{gen},i}^0(k)$ denotes the nominal (actual) set-point of the generation unit and the term $u_{\textrm{gen},i}^{\textrm{feasible}}(k)$ represents an arbitrary set-point from the set of all feasible operating points $\{\cdot\}$ to which the unit can be steered to provide operational flexibility. Both terms can be chosen to be time-variant. They are given here for time-step $k$.
The set of all feasible operation points thus depends upon the internal status of the generation unit, as defined by the terms $\xi_i(k)$, $w_i(k)$, $v_i(x_i(k))$ and $C_i x_i(k)$, and bounded by the unit's power rating constraints (Eq.~\ref{eq:general_powernode}~(b--d)).

The maximum available flexibility for up/down power regulation is given as
\begin{IEEEeqnarray}{lcl}\label{eq:FlexAssessmentPi2}
\pi^{+}_{\textrm{max},i}(k) &=& \min \left[ \eta_{\textrm{gen}} \left( \xi^{\textrm{max}} - w^{\textrm{min}} - v_x - C \dot{x} \right),\; u_{\textrm{gen}}^{\textrm{max}} \right]_{k,i} \nonumber \\ 
&-& u_{\textrm{gen},i}^{0}(k) \quad , \nonumber \\
    \pi^{-}_{\textrm{min},i}(k) &=& \max \left[ \eta_{\textrm{gen}} \left( \xi^{\textrm{min}} - w^{\textrm{max}} - v_x - C \dot{x} \right),\; u_{\textrm{gen}}^{\textrm{min}} \right]_{k,i} \nonumber \\ 
&-& u_{\textrm{gen},i}^{0}(k) \quad , 
\end{IEEEeqnarray}
in which $w^{\textrm{min}}_i(k)$ and $w^{\textrm{max}}_i(k)$ define the min./max. allowable curtailment for generation unit $i$ at time-step $k$. In case the primary fuel supply is controllable, the terms $\xi_i^{\textrm{min}}(k)$ and $\xi_i^{\textrm{max}}(k)$ define the minimum/maximum allowable primary power provision. Please note that the sign of the storage power term $C\dot{x}$ is negative when providing positive power, i.e. discharging ($C\dot{x} < 0$), and positive when providing negative power, i.e. charging ($C\dot{x} > 0$). In the time-discrete case the term $C\dot{x}$ becomes $C\delta x = C(\,x(k) - x(k-1)\,)$.

This flexibility assessment for metric $\pi$ (Eq.~\ref{eq:FlexAssessmentPi1}--\ref{eq:FlexAssessmentPi2}) can be extended to other two metrics, $\rho$ and $\epsilon$, via time-differentiation and integration, respectively. The flexibility assessment for all other power system unit types can be accomplished in a similar fashion. Please note that the maximum available flexibility calculated in this way is without any consideration of how long a certain power system unit would need to \emph{reach} a new operation point that allows this provision of flexibility. 

\subsection{Visualization of Operational Flexibility}
The three thus calculated flexibility metrics span a so-called \emph{flexibility volume}, which can be represented in its simplified form as a flexibility cube for a generic power system unit $i$, with the terms {$\pi^{+}$, $\pi^{-}$, $\rho^{+}$, $\rho^{-}$, $\epsilon^{+}$, $\epsilon^{-}$} as its vertices or extreme points. A qualitative illustration of this is shown in Fig.~\ref{fig:Flex_Cube}, where the flexibility volume is cut into eight separate sectors. 

The evolution over time of the (maximum) available operational flexibility from a generic storage unit with both load and generation terms, $u_{\textrm{load}}(k)$ and $u_{\textrm{gen}}(k)$, is illustrated in Fig.~\ref{fig:Flex_Evolution}. The plots show that the available operational flexibility is highly time-variant due to the actual storage usage over time.

\begin{figure*}[t]
%\vspace{-1.0cm}
%\begin{picture}(80,140)
%\put(0,0)
\centering               %left bottom right top
{\includegraphics[trim = 0cm   0cm   0cm   0cm, clip=true, angle=0, width=0.40\linewidth, keepaspectratio, draft=false]{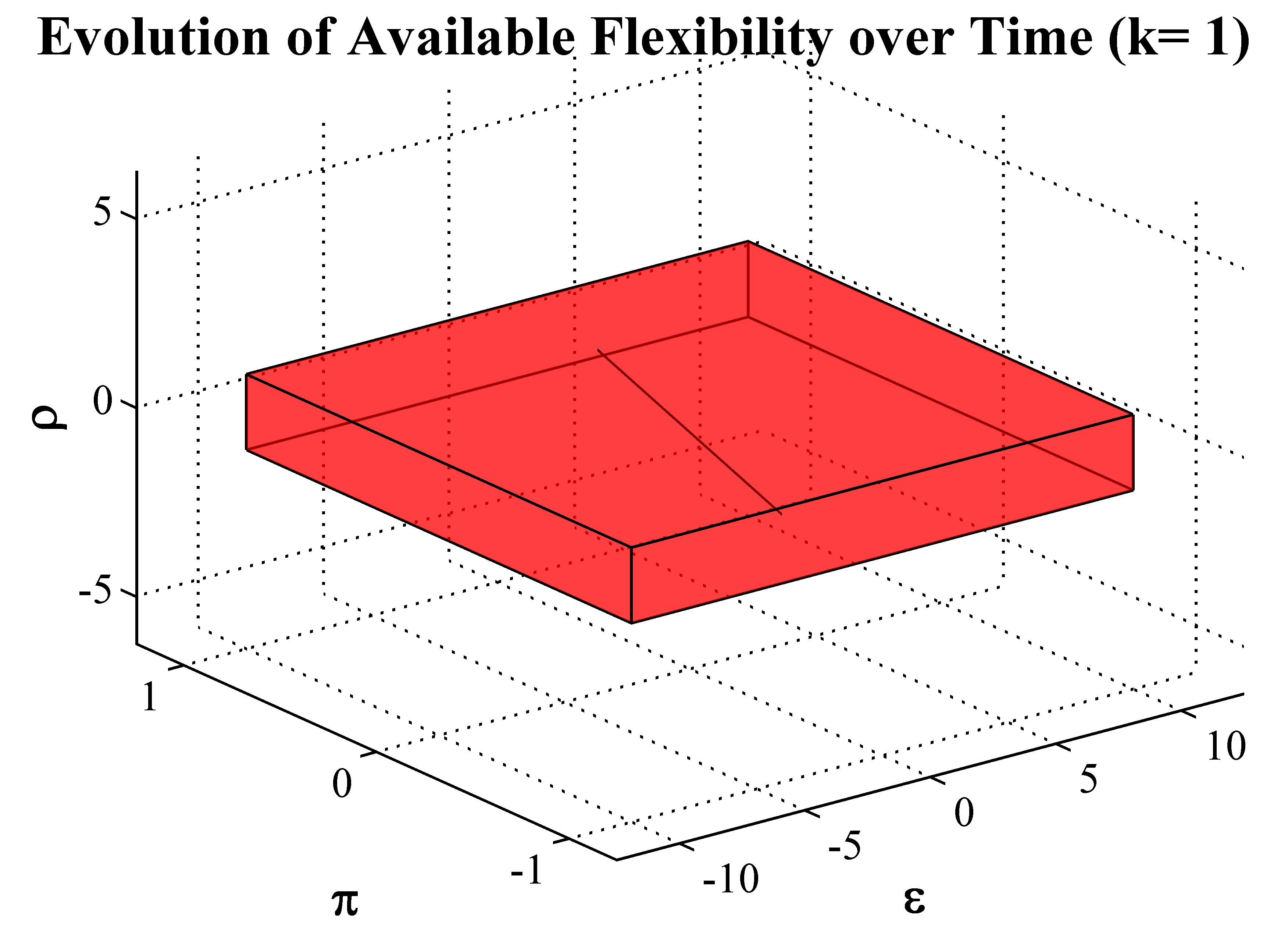}}
{\includegraphics[trim = 0cm   0cm   0cm   0cm, clip=true, angle=0, width=0.40\linewidth, keepaspectratio, draft=false]{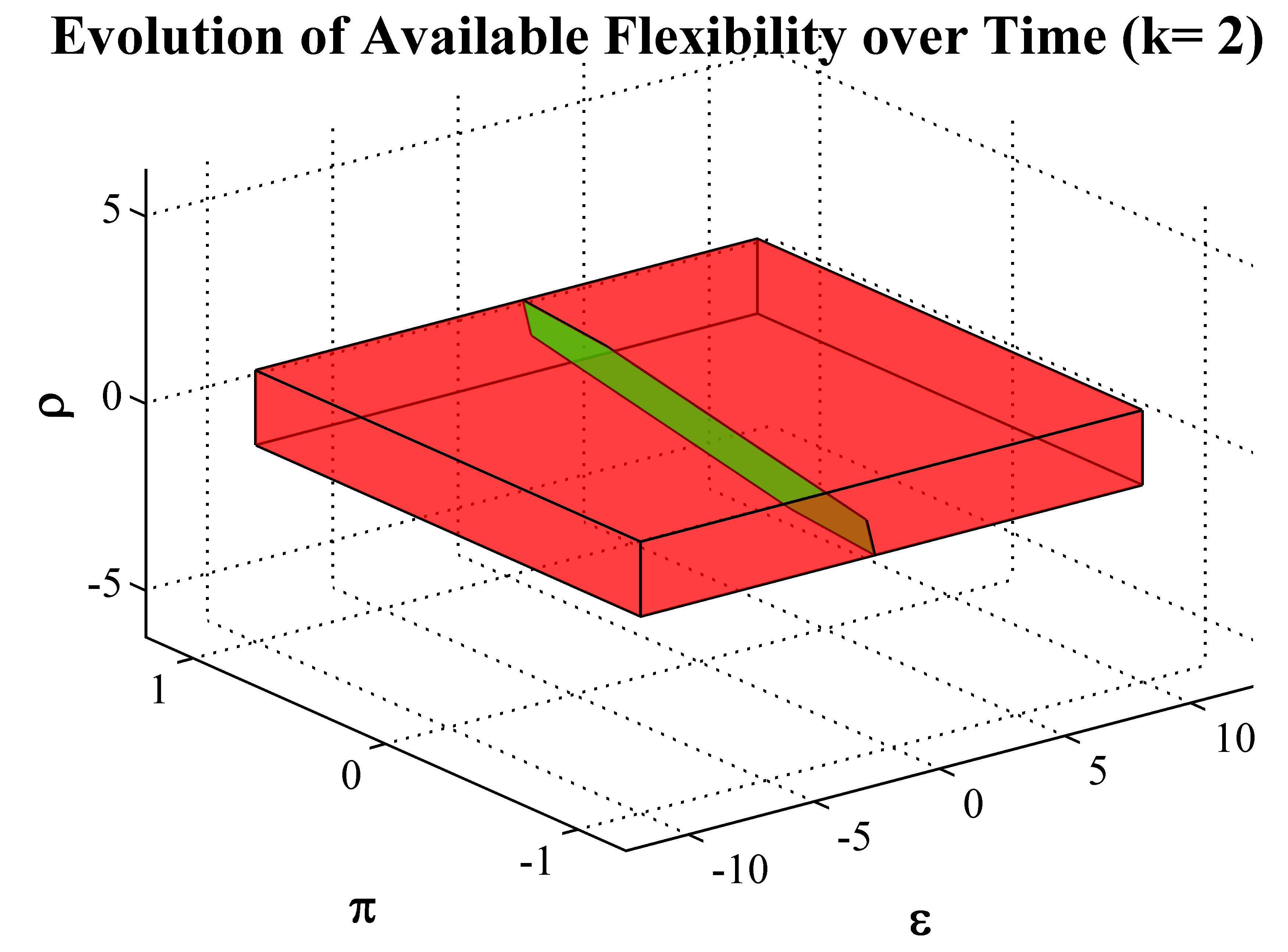}}
{\includegraphics[trim = 0cm   0cm   0cm   0cm, clip=true, angle=0, width=0.40\linewidth, keepaspectratio, draft=false]{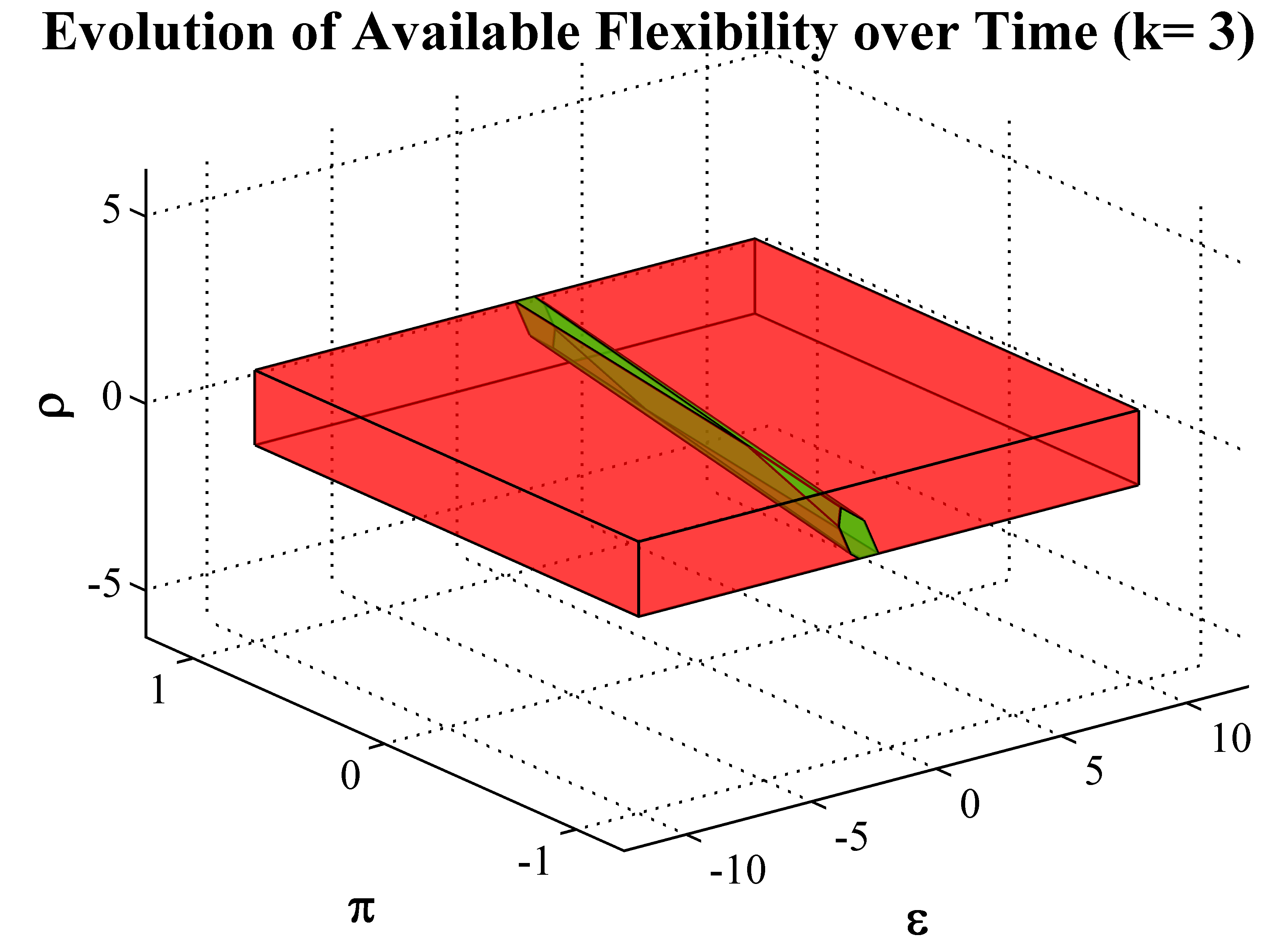}}
{\includegraphics[trim = 0cm   0cm   0cm   0cm, clip=true, angle=0, width=0.40\linewidth, keepaspectratio, draft=false]{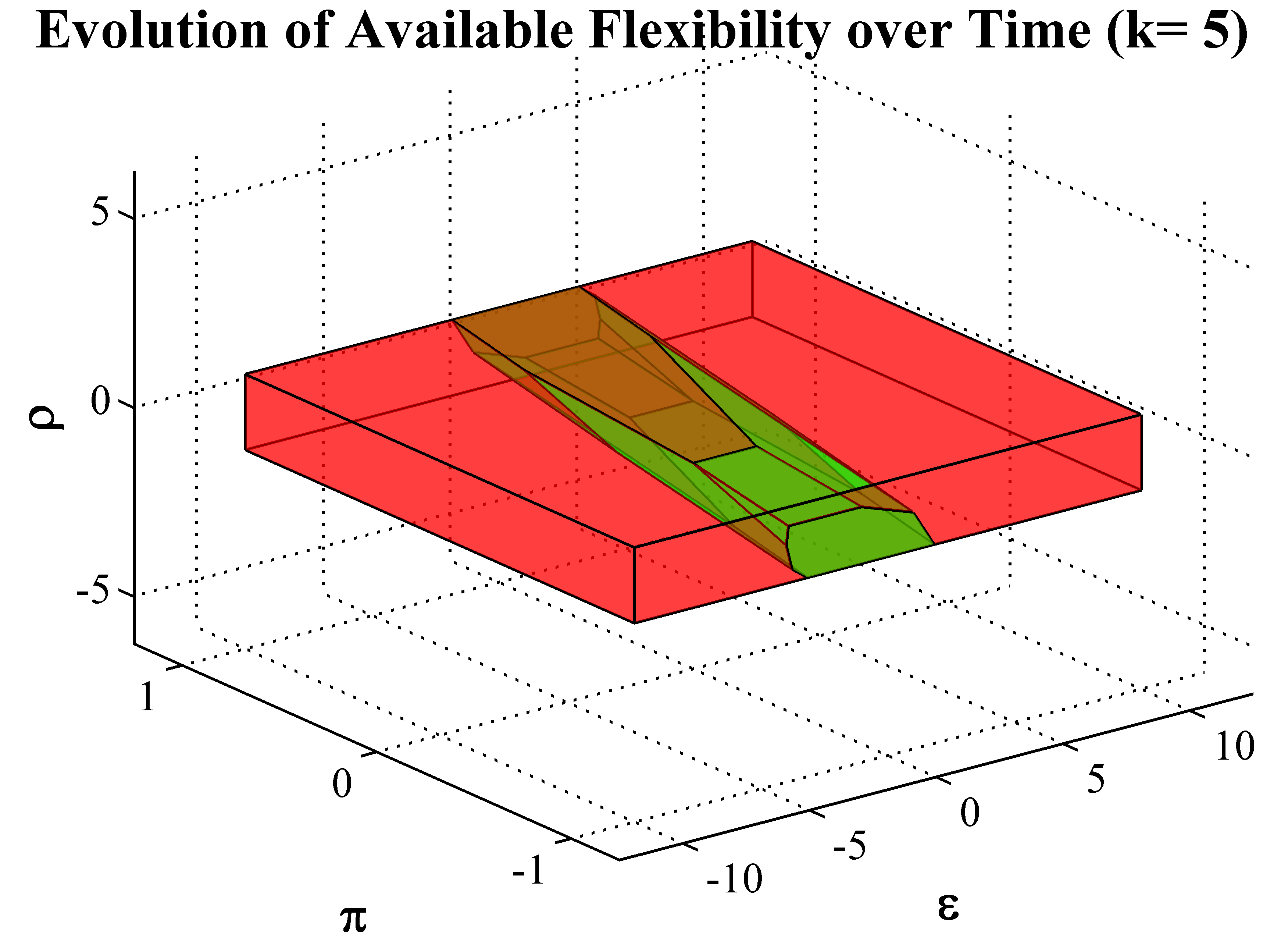}}
{\includegraphics[trim = 0cm   0cm   0cm   0cm, clip=true, angle=0, width=0.40\linewidth, keepaspectratio, draft=false]{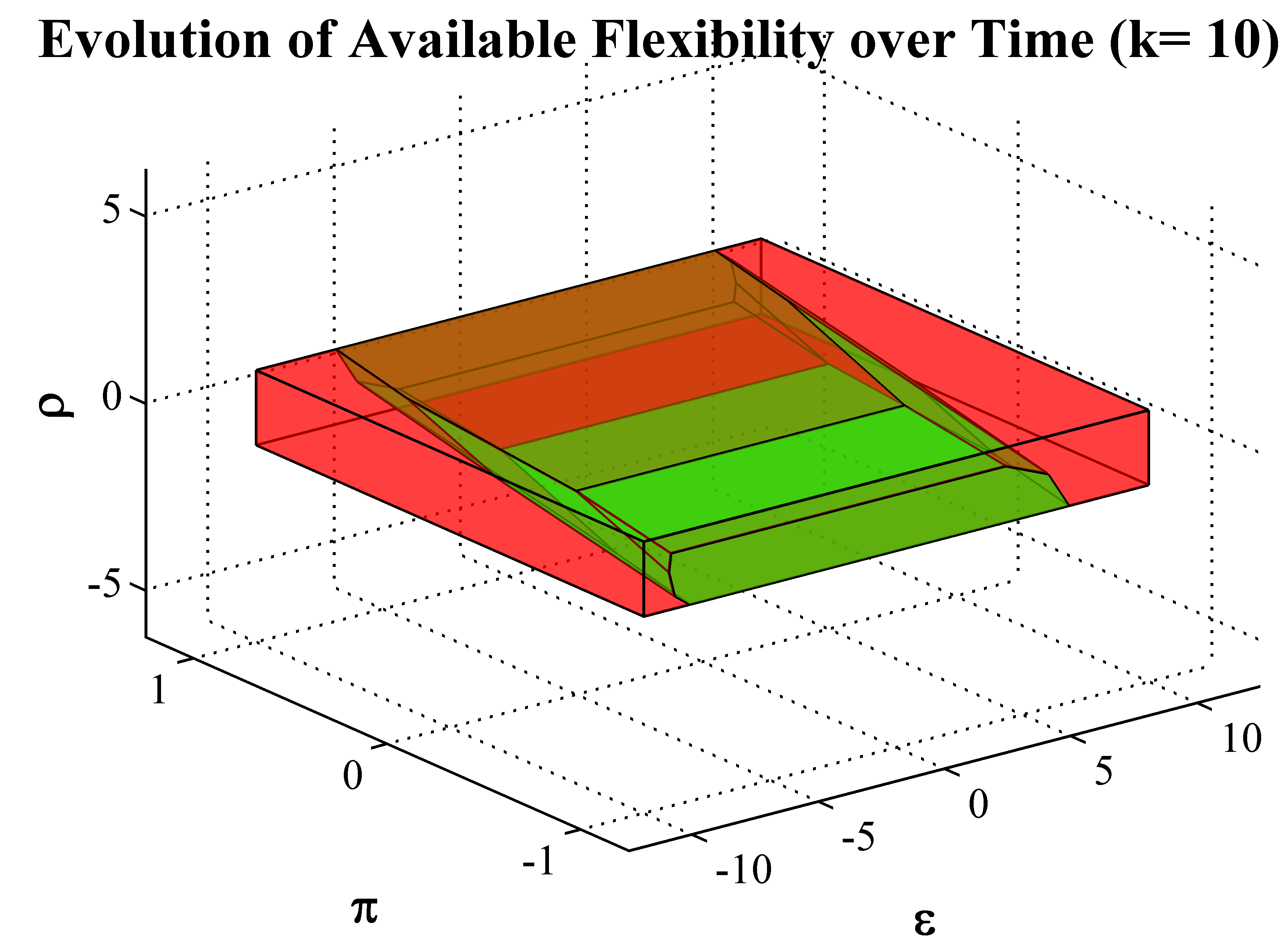}}
{\includegraphics[trim = 0cm   0cm   0cm   0cm, clip=true, angle=0, width=0.40\linewidth, keepaspectratio, draft=false]{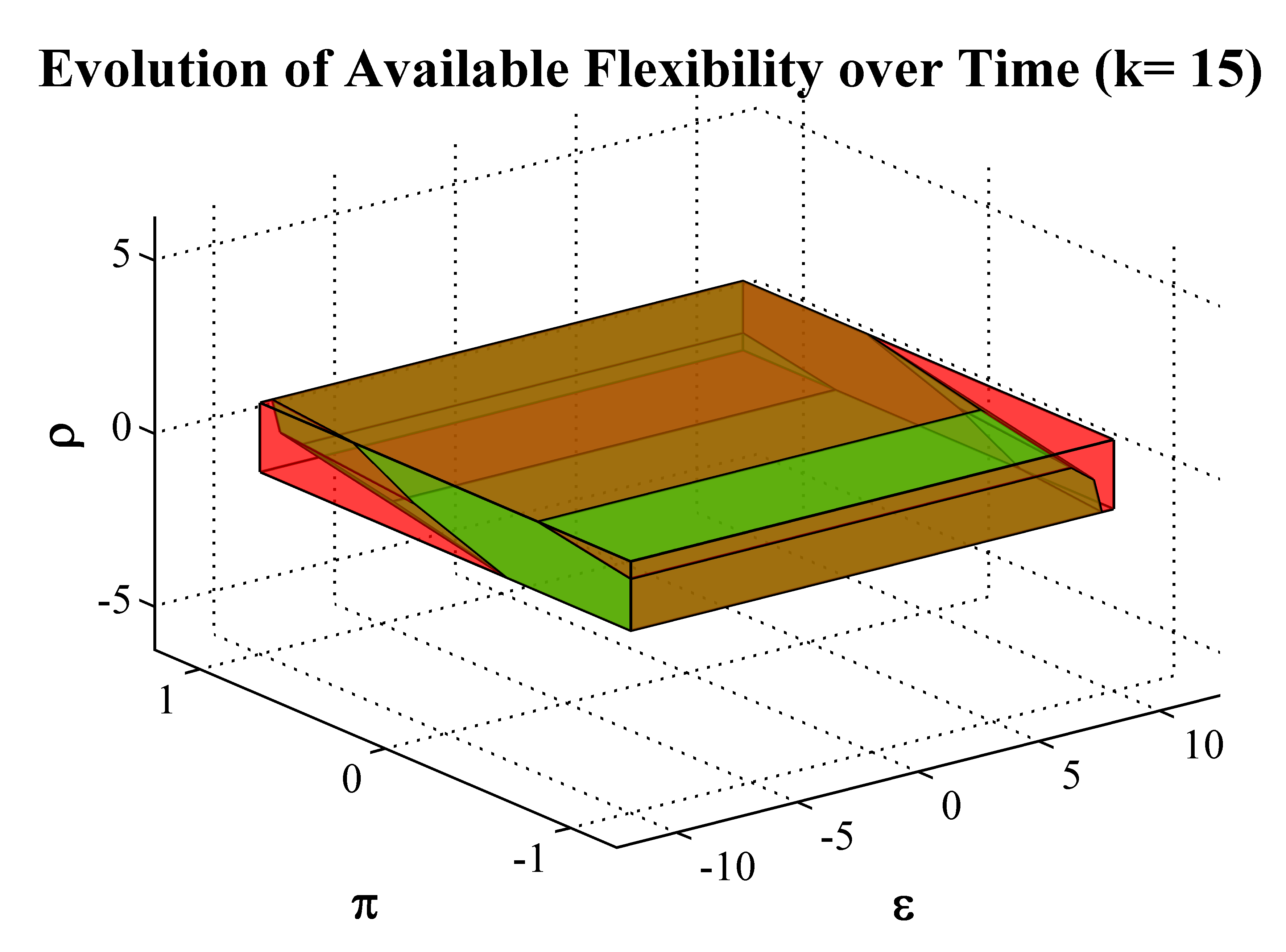}}
%\end{picture}
\caption{Evolution of available operational flexibility from a storage unit at its planned operation point ($k=0$). \newline 
Green: Time-evolution of available flexibility after k = 1h, 2h, 3h, 5h, 10h, 15h (calculated via reach sets). \newline Red: Maximum available operational flexibility at $k\rightarrow\infty$ (calculated using Eq.~(\ref{eq:FlexAssessmentPi1}--\ref{eq:FlexAssessmentPi2})).} \label{fig:Flex_Reach}
\vspace{-0.4cm}
\end{figure*}

However, when taking into account the internal double-integrator dynamics, the flexibility volume becomes a significantly more complex polytope object. An illustration of this more realistic polytope flexibility volume is given by Fig.~\ref{fig:Flex_Reach}. Here, the information of how long it takes to reach a certain new operation point providing a required set $\{\rho,\pi,\epsilon\}$ of operational flexibility is explicitly given. The set of reachable operation points providing additional flexibility (green) becomes larger when the available time span is longer. The flexibility set converges towards the set of maximal flexibility (red) as defined by the underlying technical constraints of a given power system unit. Calculating the available set of operational flexibility that is achievable after a given number of time-steps $k$ is equivalent to a classical \emph{reach set} calculation. This later approach, although more exact, is significantly more computationally expensive than the simpler analytic approach sketched out previously by Eq.~(\ref{eq:FlexAssessmentPi1}--\ref{eq:FlexAssessmentPi2}). 

%For our calculations, we have used the reachability functions of the \texttt{MPT Toolbox}~\cite{mpt} in~\texttt{Matlab}. %, which involves besides other things the calculation of the \emph{controllability gramian} $W_C$.
For the reach-set calculations, the reachability functions of the \texttt{MPT Toolbox}~\cite{mpt} have been used in~\texttt{Matlab}. There a so-called polytope method is employed that involves besides other things the calculation of the \emph{Controllability Gramian} $W_C$. (See~\cite[p.~19~ff.]{ET_manual} for a general discussion of reachability analysis.) 
The advantage of the \texttt{MPT Toolbox} is that it explicitly allows the usage of box constraints for inputs and states of dynamical systems. In power systems, a typical example of a box constraint are the limitations on min/max power ramp-rate, e.g.~$\dot{u}_{\textrm{gen.}}^{\textrm{min}} \le \dot{u}_{\textrm{gen.}}\,(k) \le \dot{u}_{\textrm{gen.}}^{\textrm{max}}$, and min/max power output, e.g.~$u_{\textrm{gen.}}^{\textrm{min}} \le u_{\textrm{gen.}}\,(k) \le u_{\textrm{gen.}}^{\textrm{max}}$. 
%\begin{equation}\label{equ:ReachFlex}
% ...
%\end{equation}
%
Other approaches for calculating gramians and the corresponding reach-sets include \ac{lmi} methods, as explained in~\cite{boyd1994LMI}, as well as so-called ellipsoidal methods, which have been implemented for example in the \texttt{Ellipsoidal Toolbox}~\cite{ET}.
%Please note that ellipsoidal methods have a potential disadvantage as they only allow ellipsoidal constraints on system input and states. \\
%On the other hand they are computationally much less expensive than \texttt{MPT}'s polytope method when it comes to larger system sizes as is nicely illustrated in~\cite[p.~63~ff.]{ET_manual}.

Please note that the theoretical maximum reachability volume calculated by the  analytic approach may in fact never be fully reached by the power system unit, when using the reach set approach (Fig.~\ref{fig:Flex_Reach}). %--\ref{fig:Flex_ReachVolume}). 
This gap is due mainly to the non-infinite discrete sampling time in combination with somewhat pathological operation points at some of the flexiblity cube's vertices, e.g.~fully discharging a storage unit ($\pi^{-}$) while at the same time keeping it at its maximum energy storage level ($\epsilon^{+}$). 
%This gap is due mainly to the non-infinite discrete sampling time in combination with somewhat pathological operation points at some of the flexiblity cube's vertices. These operation points exist but are not of any practical relevance, e.g.~fully discharging a storage unit ($\pi^{-}$) while at the same time keeping it at its maximum energy storage level ($\epsilon^{+}$). 
%\begin{equation}\label{equ:ReachFlex}
% ...
%\end{equation}

%For the sake of simplicity in representing flexibility volumes we will stick to the simplified flexibility cubes for the remainder of this paper. 

\subsection{Aggregation of Operational Flexibility}

\begin{figure*}[t]
%\vspace{-1.0cm}
%\begin{picture}(80,140)
%\put(0,0)
\centering               %left bottom right top
{\includegraphics[trim = 0.25cm   0cm   0.0cm   0cm, clip=true, angle=0, width=0.325\linewidth, keepaspectratio, draft=false]{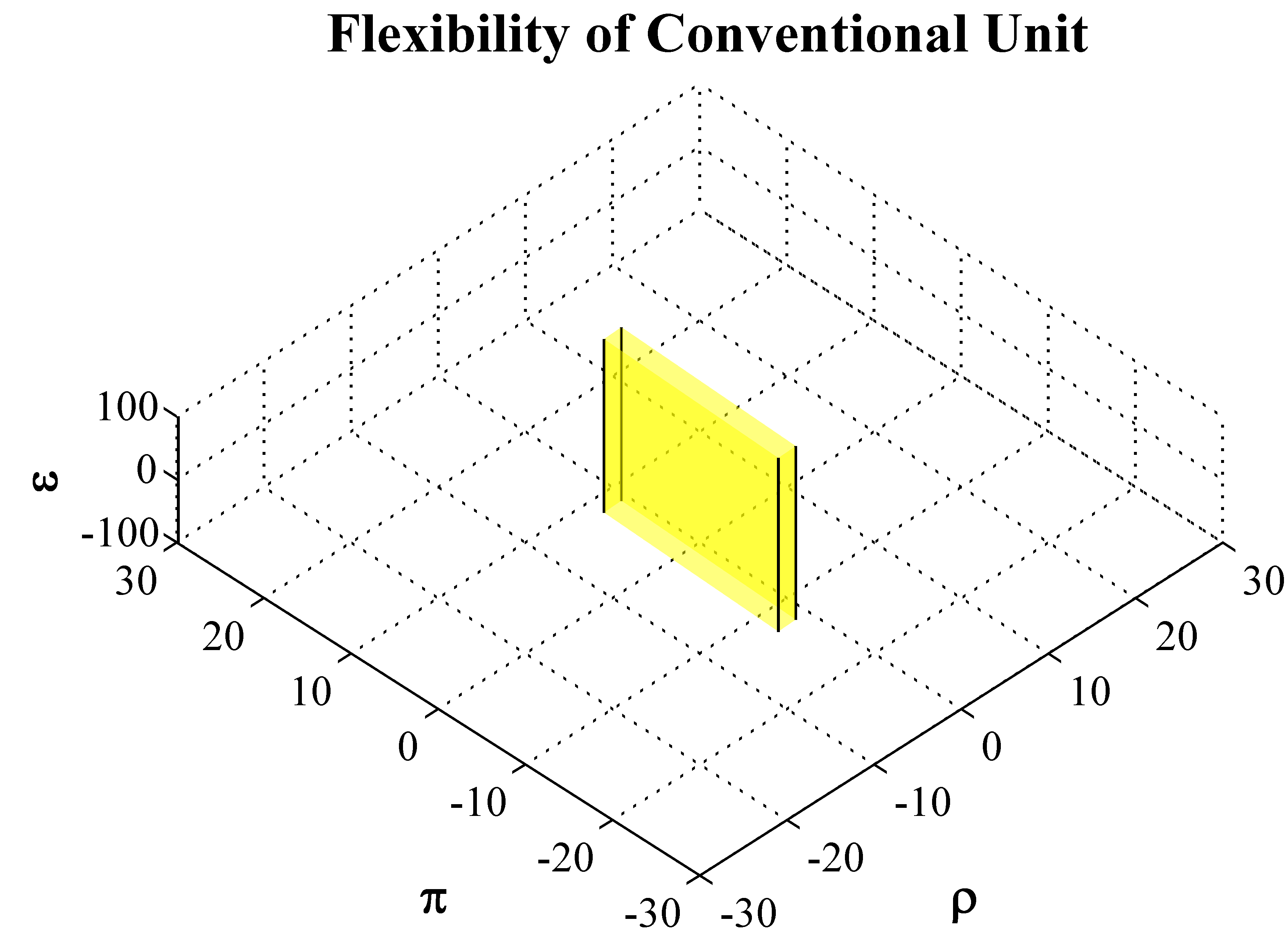}}
{\includegraphics[trim = 0.25cm   0cm   0.0cm   0cm, clip=true, angle=0, width=0.325\linewidth,  keepaspectratio, draft=false]{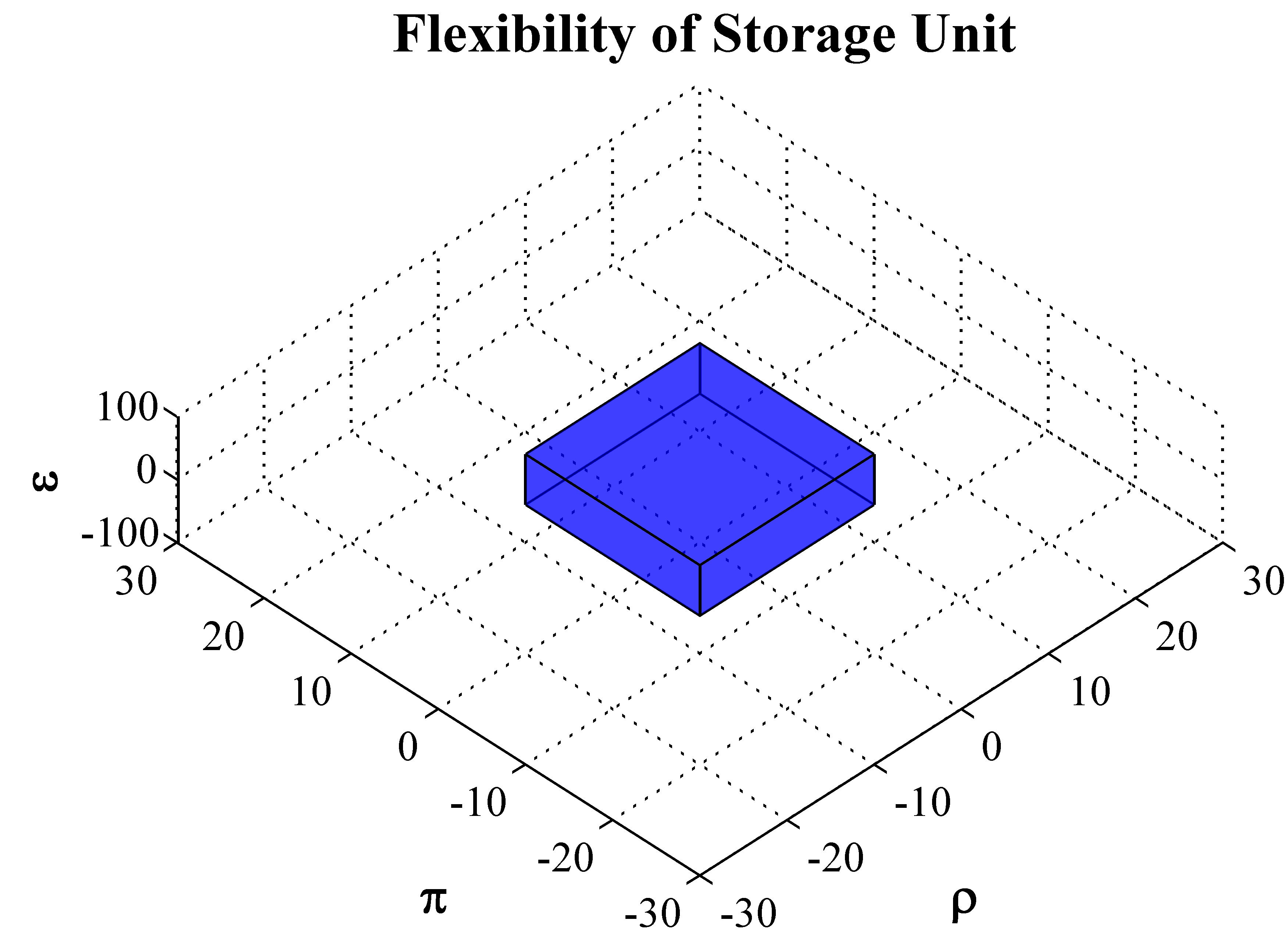}}
{\includegraphics[trim = 0.25cm   0cm   0.0cm   0cm, clip=true, angle=0, width=0.325\linewidth,  keepaspectratio, draft=false]{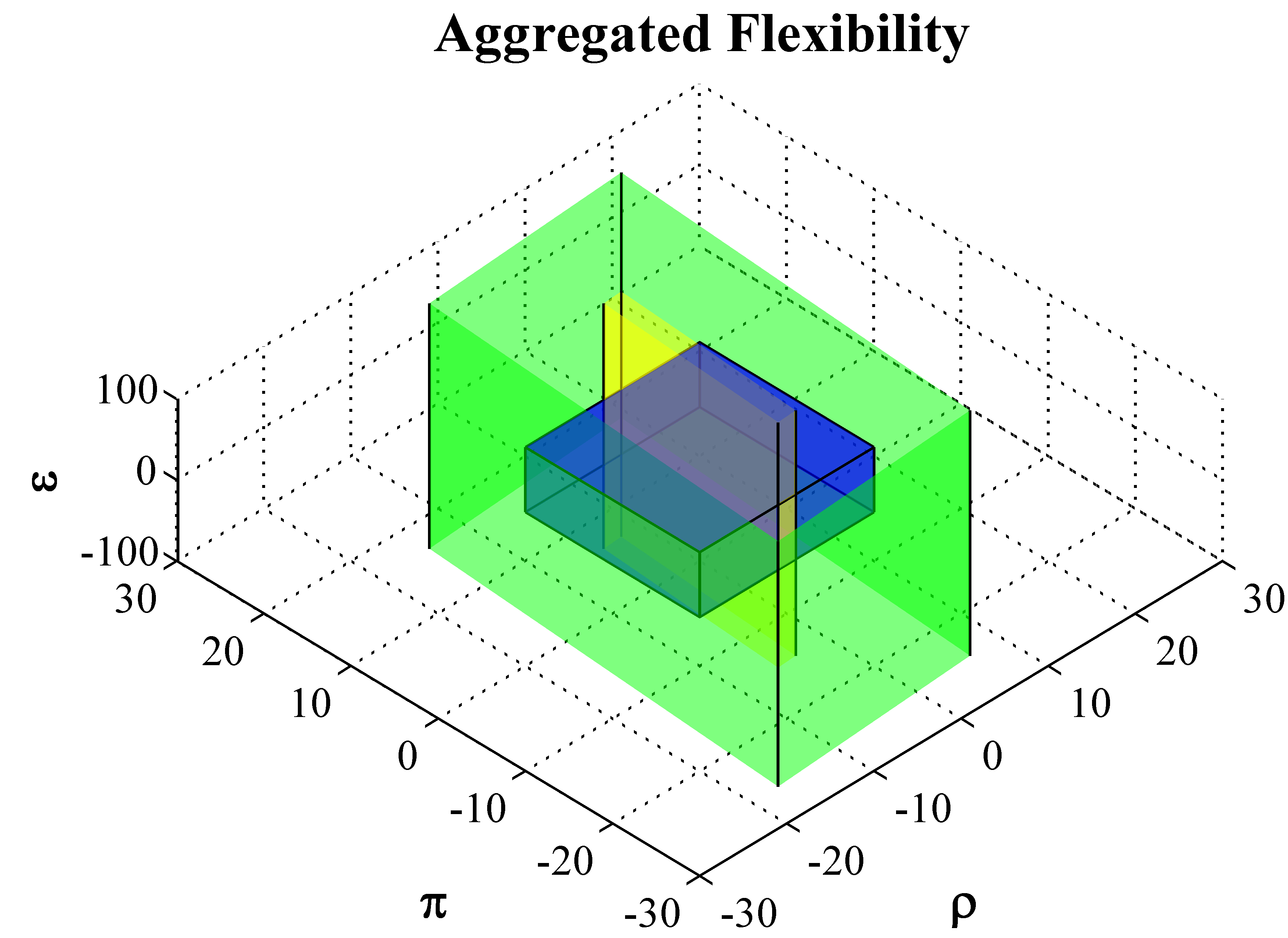}}
\\
\vspace{0.5cm}
{\includegraphics[trim = 0.25cm   0cm   0.0cm   0cm, clip=true, angle=0, width=0.325\linewidth, keepaspectratio, draft=false]{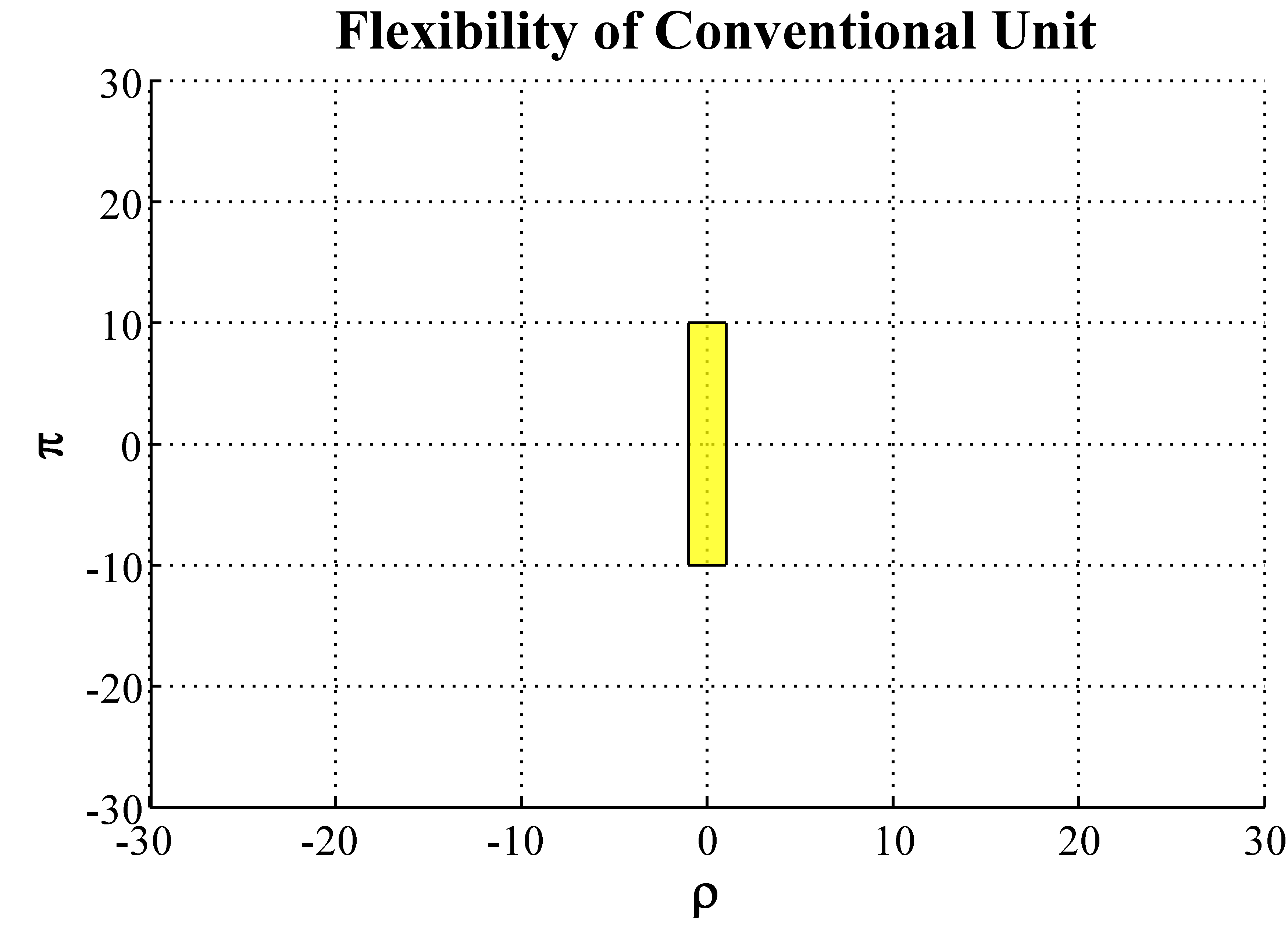}}
{\includegraphics[trim = 0.25cm   0cm   0.0cm   0cm, clip=true, angle=0, width=0.325\linewidth,  keepaspectratio, draft=false]{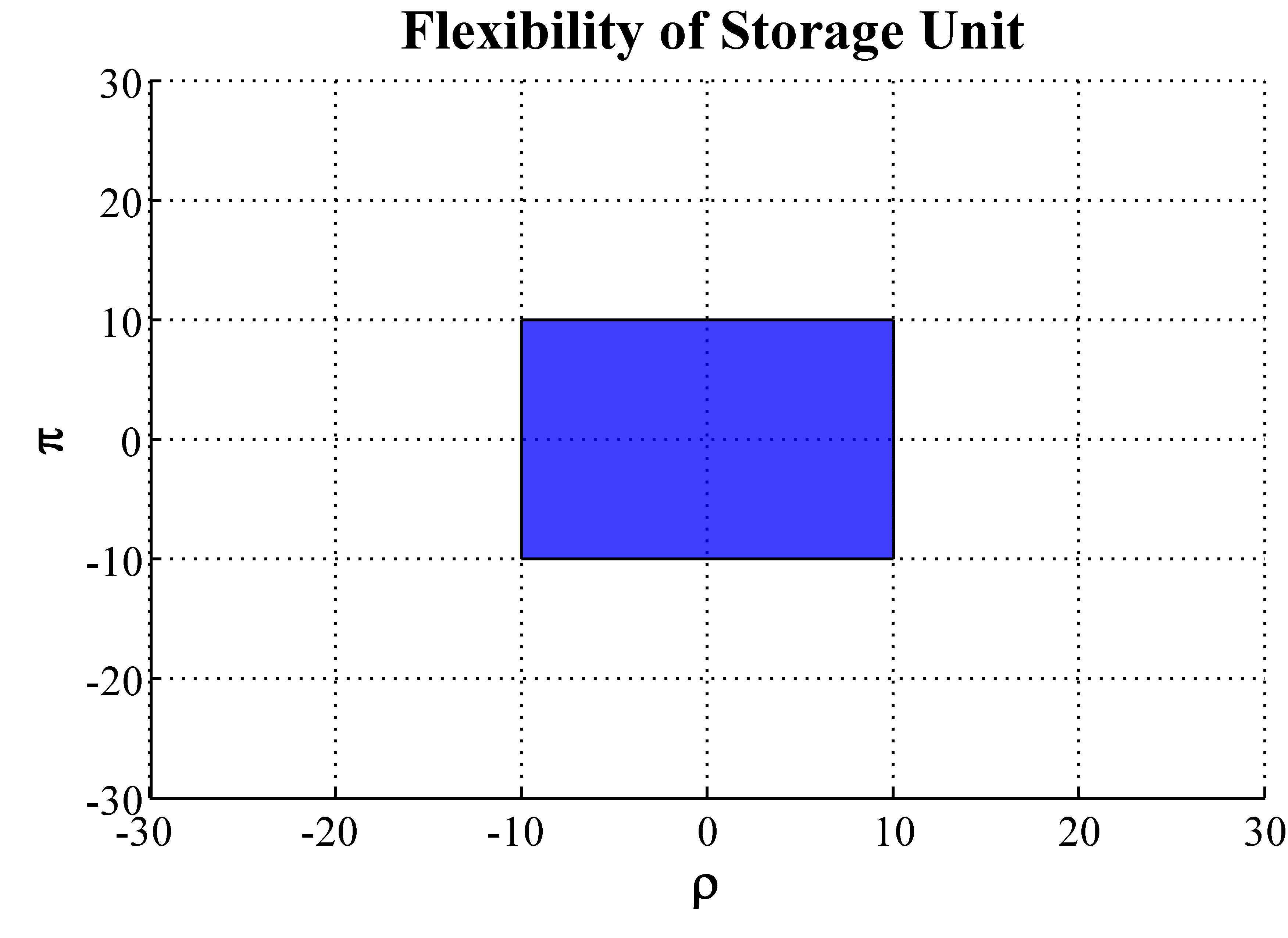}}
{\includegraphics[trim = 0.25cm   0cm   0.0cm   0cm, clip=true, angle=0, width=0.325\linewidth,  keepaspectratio, draft=false]{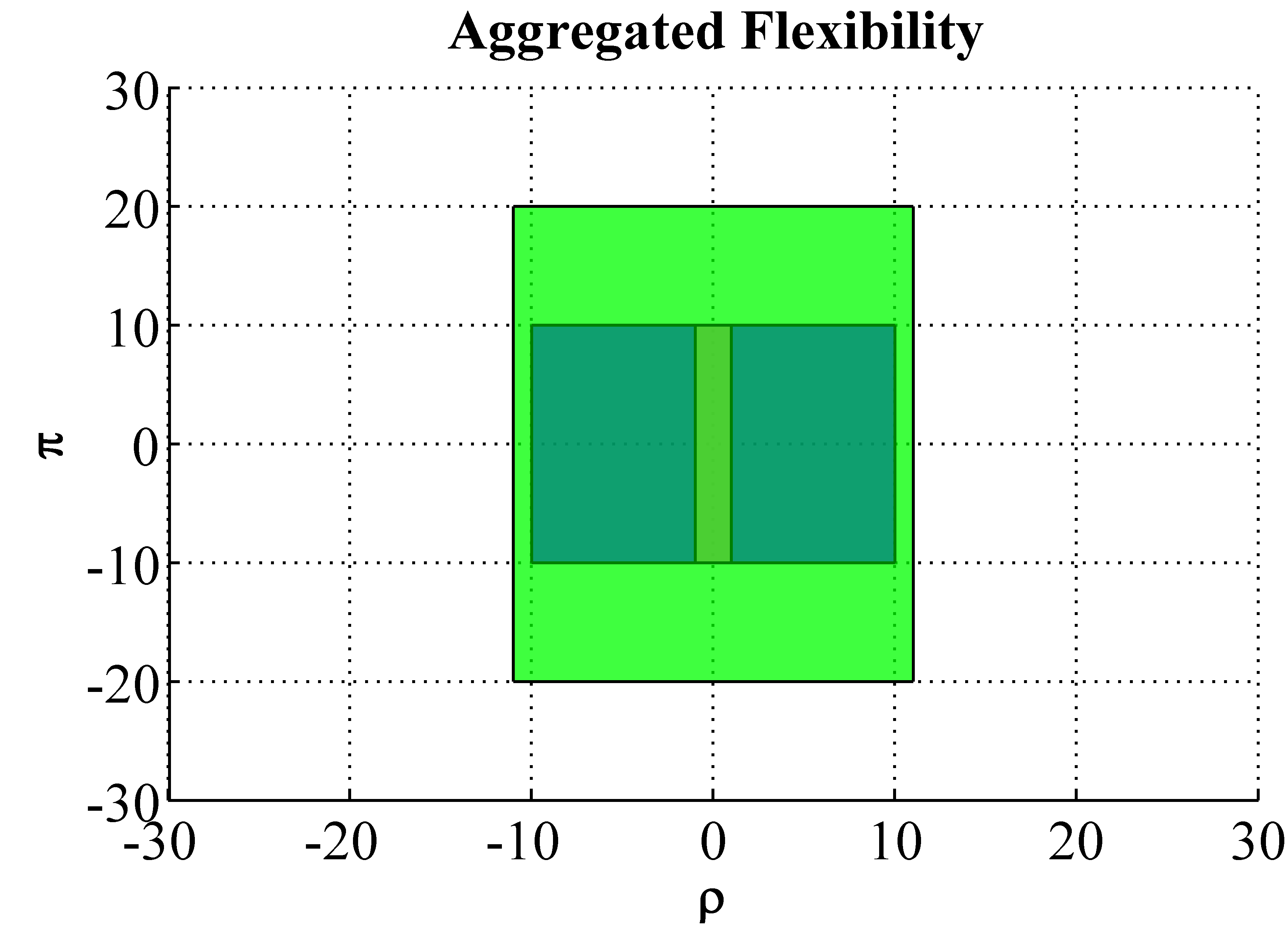}}
%\end{picture}
\caption{Aggregation of maximum operational flexibility of individual power system units.\newline Flexibility of conventional unit with no energy constraint (yellow), flexibility of energy-constrained storage (blue) and aggregated flexibility of both units (green).} \label{fig:Flex_Aggregation}
%\vspace{-0.5cm}
\end{figure*}

An important question in power system analysis is how a group or pool of power system units act together in achieving a given objective, i.e.~delivering a scheduled power trajectory or providing ancillary services by tracking a control signal. Pooling together different power system units to provide a service that they cannot provide individually is an active research field. A prime example is to combine a dynamically slow power plant with a dynamically fast, but energy-constrained storage unit to provide fast frequency regulation that neither of the units could provide individually~\cite{JinLuLu:2013} due to the lack of one flexibility metric, i.e.~the missing fast ramp-rate capability $\rho$ of the power plant, or another, i.e.~the small energy capability $\epsilon$ of the storage unit.
Obtaining the aggregated operational flexibility that a pool of different power system units provides, is equivalent to aggregating the flexibility volumes of the individual units. Since these are given by more or less complex polytope sets, depending on the chosen calculation approach presented in the previous section, a well-known polytope operation, the \emph{Minkowski Summation}, can be employed for calculating the aggregated flexibility of the pool. In the following, we illustrate the aggregation of a slow-ramping power plant together with a fast-ramping but energy-constrained storage unit in Fig.~\ref{fig:Flex_Aggregation}. We assume that within the grid zone of a unit pool, grid constraints are minor and not of practical relevance for the quantification of aggregated flexibility. Although this is a simplifying assumption, it is often used, e.g.~in power markets operation.

The aggregation of  two or more power system units leads to the addition of individual flexibility metrics
\begin{equation}
\lbrace \rho, \pi, \epsilon \rbrace_{\textrm{agg}} = \lbrace \rho, \pi, \epsilon \rbrace_{\textrm{slow}} + \lbrace \rho, \pi, \epsilon \rbrace_{\textrm{fast}} \quad .
\end{equation}
The aggregation of the operational flexibility of both units, given individually by their respective polytope objects, is accomplished via \emph{Minkowski Summation}
%\begin{eqnarray}
%	\rho_{\textrm{agg}}^{+} &= \sum\limits_i \rho_{i}^{+}  \; , \ \	\rho_{\textrm{agg}}^{-} &= \sum\limits_i \rho_{i}^{-} \nonumber \\
%	\pi_{\textrm{agg}}^{+} &= \sum\limits_i \pi_{i}^{+} \; , \ 		\pi_{\textrm{agg}}^{-} &= \sum\limits_i \pi_{i}^{-} \\
%	\epsilon_{\textrm{agg}}^{+} &= \sum\limits_i \epsilon_{i}^{+} \; , \ \	
%	\epsilon_{\textrm{agg}}^{-} &= \sum\limits_i \epsilon_{i}^{-} \nonumber \quad . \nonumber
%\end{eqnarray}
\begin{eqnarray}
	\rho_{\textrm{agg}}^{+} &=& \sum\limits_i \rho_{i}^{+}  \; , \quad	
	\rho_{\textrm{agg}}^{-} \ \ = \ \ \sum\limits_i \rho_{i}^{-} \quad , \nonumber \\
	\pi_{\textrm{agg}}^{+}  &=& \sum\limits_i \pi_{i}^{+} \; , \quad 		
	\pi_{\textrm{agg}}^{-}  \ \ = \ \ \sum\limits_i \pi_{i}^{-} \quad , \\
	\epsilon_{\textrm{agg}}^{+} &=& \sum\limits_i \epsilon_{i}^{+} \; , \quad \;	
	\epsilon_{\textrm{agg}}^{-} \ \ = \ \ \sum\limits_i \epsilon_{i}^{-} \quad\, . \nonumber
\end{eqnarray}
The slow-ramping unit, e.g.~a thermal power plant, with $\lbrace \rho, \pi, \epsilon \rbrace_\textrm{slow}$, is assumed to have an unlimited fuel supply, which implies that no energy constraints exist and that the energy provision capability is infinite ($\epsilon_\textrm{slow} = \infty$). Also, the potential power output $\pi$ is large. Dynamically slow means in this context that the power ramp-rate $\rho$ is small. The fast-ramping storage unit, e.g.~a fly-wheel or battery system, with $\lbrace \rho, \pi, \epsilon \rbrace_\textrm{fast}$, has a limited run-time bounded by energy constraints of the storage unit and thus only a limited energy storage capability exists ($0 < \epsilon_\textrm{fast} \ll \infty$). 

As is often the case for storage units, ramp-rate $\rho$ is large whereas power capability $\pi$ is comparatively small. Depending on storage technology, time-dependent storage losses, $v(x)$, can be significant. This is notably the case of fly-wheel energy storage systems, where storage losses due to bearing friction become large when going beyond a storage cycle duration of a few minutes.

\subsection{Available versus Needed Operational Flexibility}

At last we compare the \emph{needed} operational flexibility for mitigating a disturbance event with the \emph{available} flexibility that a given power system can offer.
The needed flexibility could, for example, be derived from the assumed worst-case succession of the combined wind and \ac{pv} in-feed forecast errors over a given time interval (see~\cite{Makarov:2009} as an illustration of needed flexibility for balancing wind forecast errors in the CAISO grid).

\begin{figure*}[t]
%\vspace{-1.0cm}
%\begin{picture}(80,140)
%\put(0,0)
\centering               %left bottom right top
{\includegraphics[trim = 0.25cm   0cm   0cm   0cm, clip=true, angle=0, width=0.325\linewidth, keepaspectratio, draft=false]{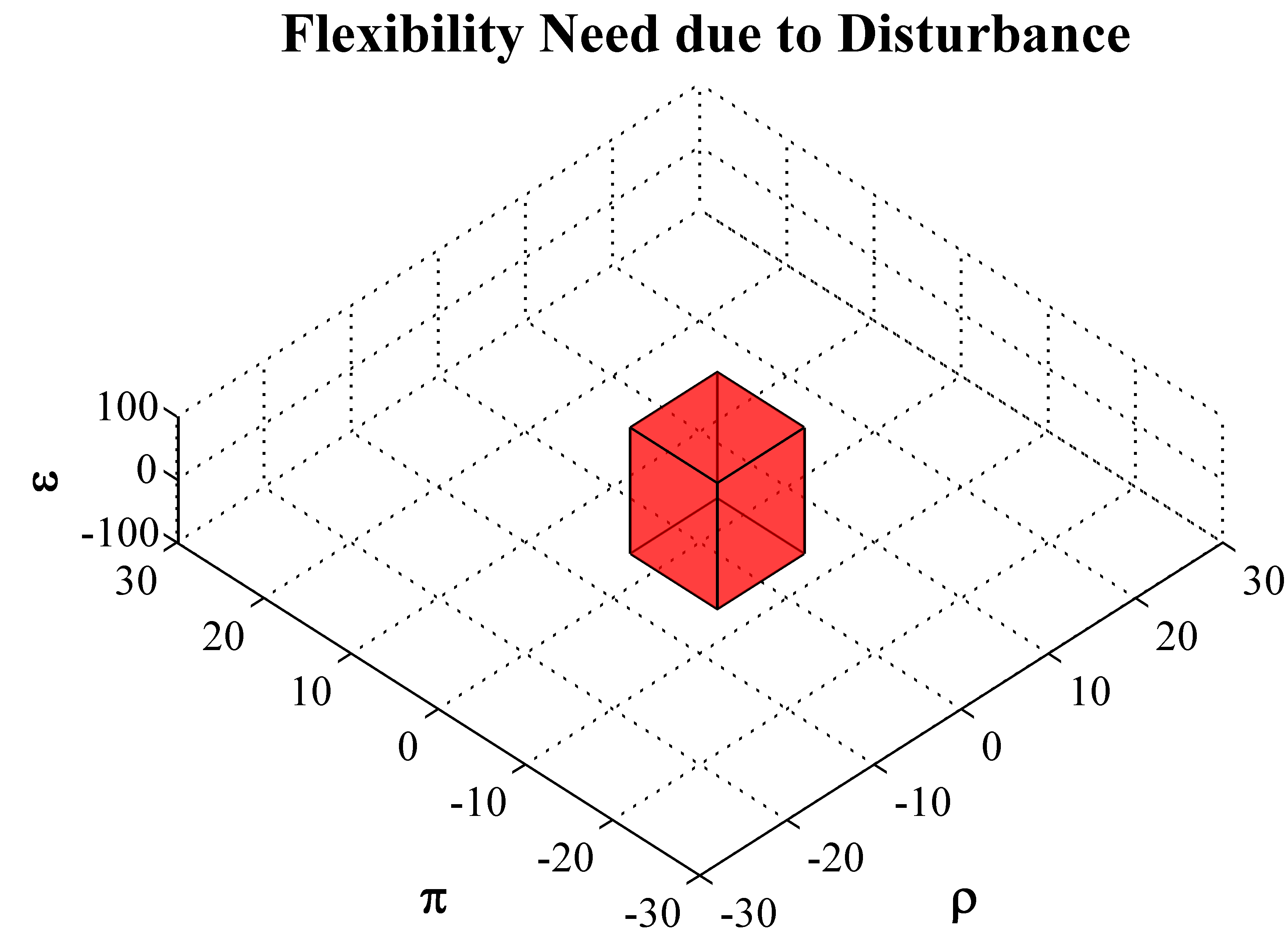}}
{\includegraphics[trim = 0.25cm   0cm   0cm   0cm, clip=true, angle=0, width=0.325\linewidth, height=0.243\linewidth, draft=false]{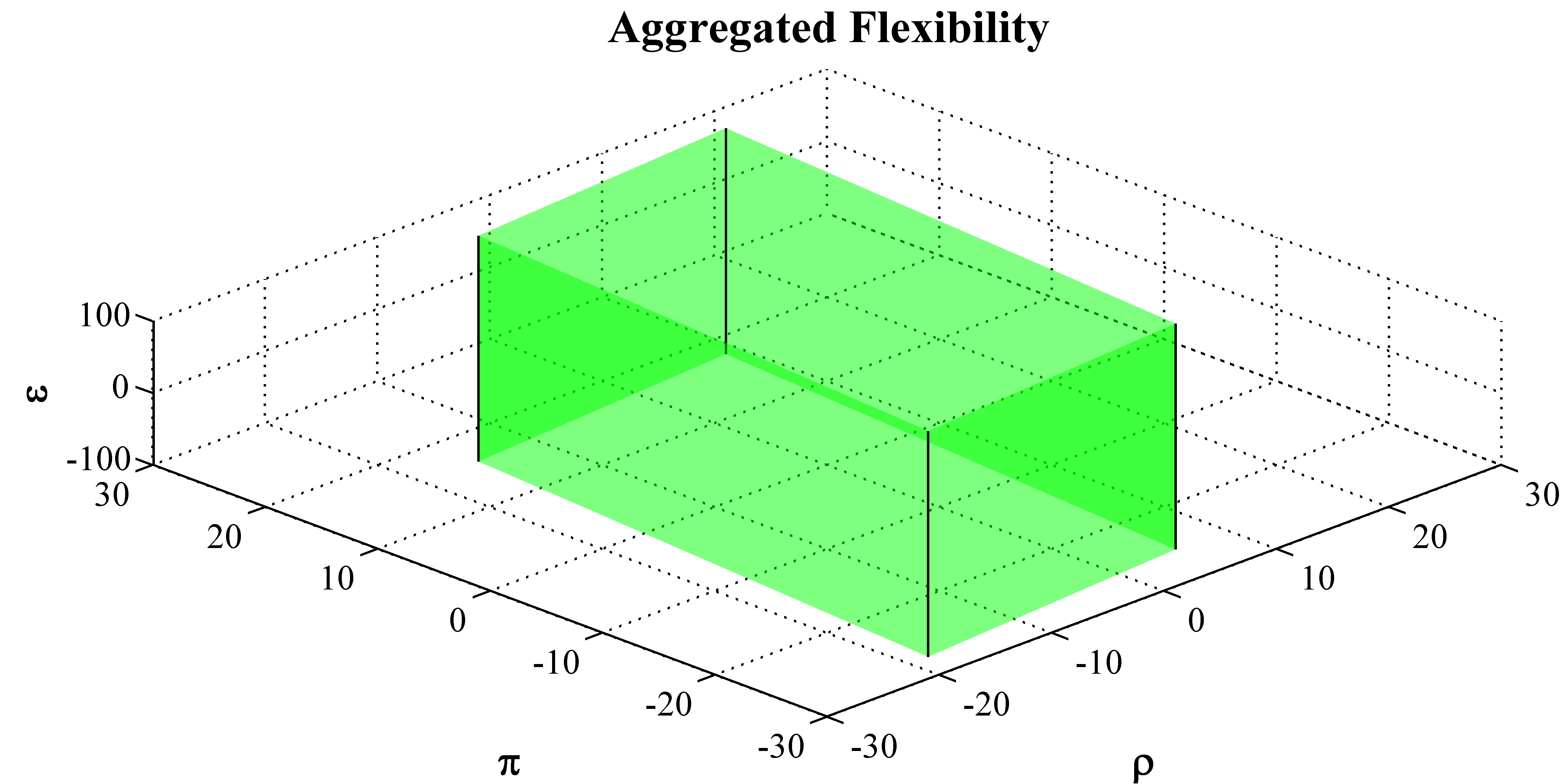}}
{\includegraphics[trim = 0.25cm   0cm   0cm   0cm, clip=true, angle=0, width=0.325\linewidth,  keepaspectratio, draft=false]{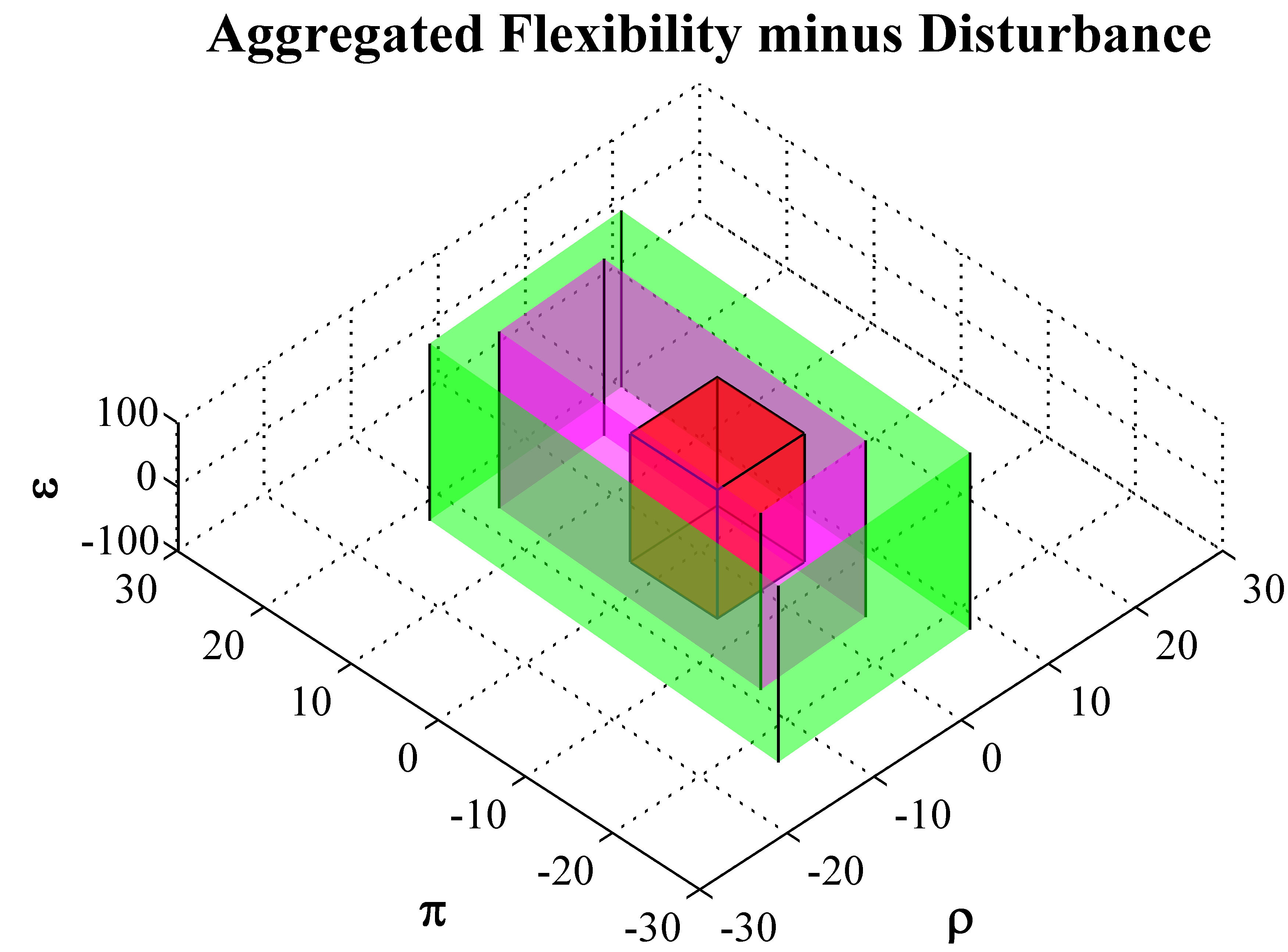}}
\\
\vspace{0.5cm}
{\includegraphics[trim = 0.25cm   0cm   0cm   0cm, clip=true, angle=0, width=0.325\linewidth, keepaspectratio, draft=false]{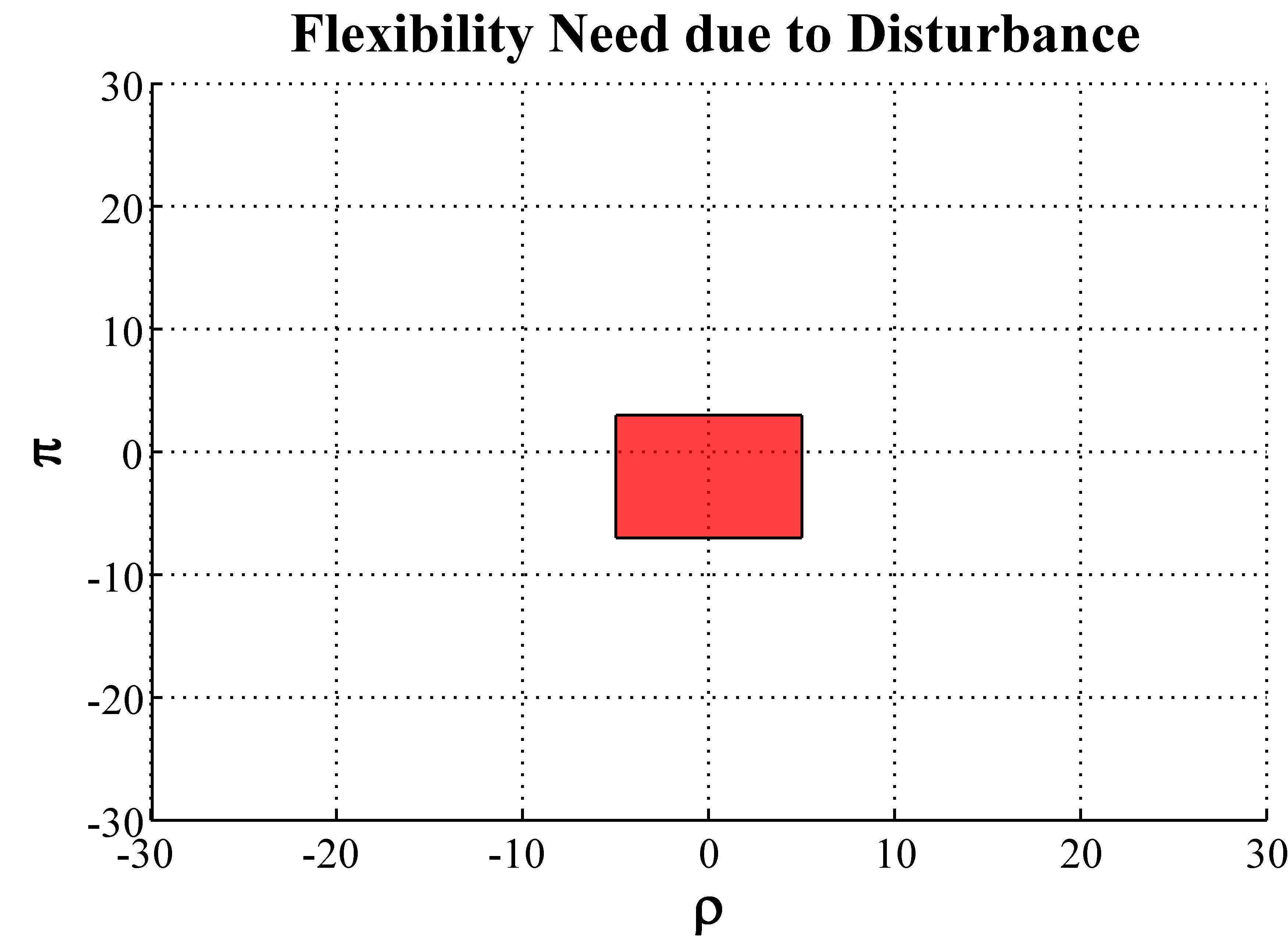}}
{\includegraphics[trim = 0.25cm   0cm   0cm   0cm, clip=true, angle=0, width=0.325\linewidth, height=0.243\linewidth, draft=false]{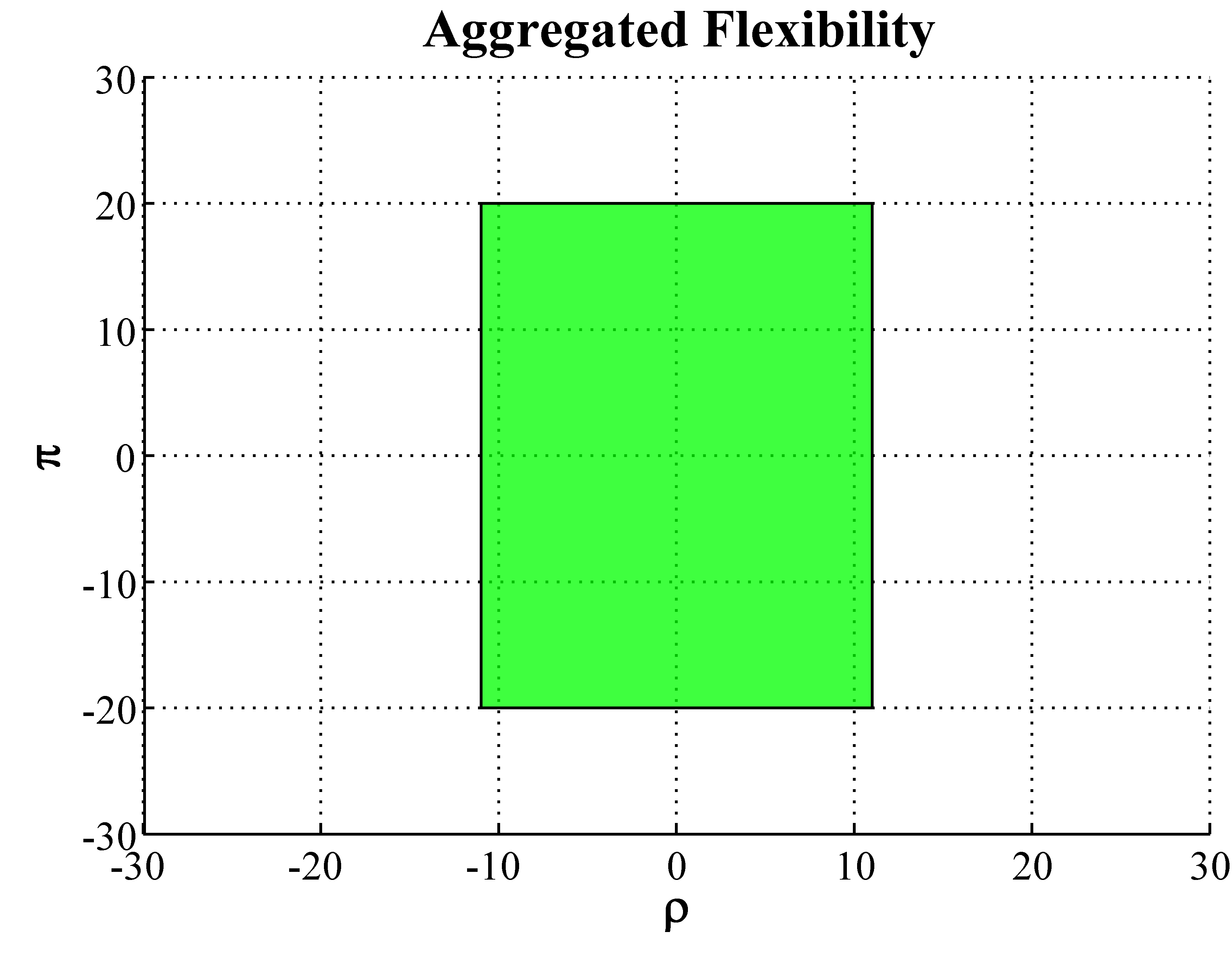}}
{\includegraphics[trim = 0.25cm   0cm   0cm   0cm, clip=true, angle=0, width=0.325\linewidth,  keepaspectratio, draft=false]{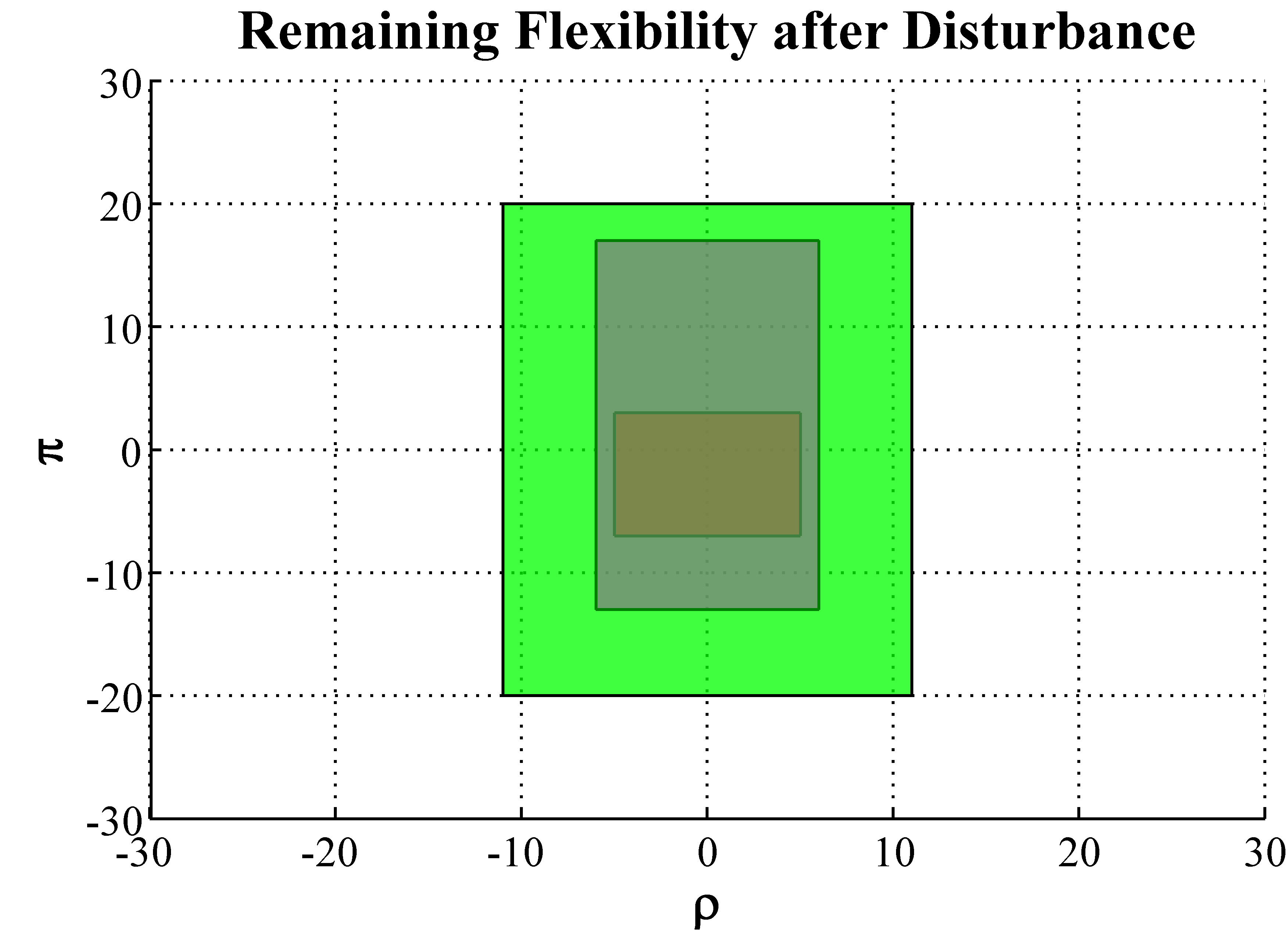}}
%\end{picture}
\caption{\emph{Needed} operational flexibility versus \emph{available} operation flexibility. \newline 
Needed flexibility volume for balancing a disturbance (red), available flexibility volume (green) and remaining flexibility volume after subtracting the needed flexibility volume (magenta).} \label{fig:Flex_Subtraction}
%\vspace{-0.5cm}
\end{figure*}

Effectively balancing this requires the ability to follow steep power ramps as well as to provide large amounts of regulating power and energy over time.
In order for a given power system to successfully accommodate such a disturbance event, the \emph{available} flexibility volume should always envelope the \emph{needed} flexibility volume, as shown in Fig.~\ref{fig:Flex_Subtraction}. If this would not be the case, flexibility capability is lacking along at least one axis of the flexibility metrics (e.g.~$\pi^{+}_{\textrm{agg.}}=0$). The disturbance event could not be fully accommodated. Calculating the polytope of the still available operational flexibility that remains while mitigating the expected disturbance boils down to another polytope operation, the \emph{Pontryagin Difference}. 

%\begin{figure}[t]
%%\vspace{-1.0cm}
%%\begin{picture}(80,140)
%%\put(0,0)
%\centering               %left bottom right top
%{\includegraphics[trim = 0cm   0cm   0cm   0cm, clip=true, angle=0, width=1.00\linewidth, keepaspectratio, draft=false]{Volume_Evolution_HQ.png}}
%%\end{picture}
%\caption{Reached flexibility volume, $V(k) = ( \rho \cdot \pi \cdot \epsilon )(k)$, after $k$ time-steps for a storage unit at its planned operation point ($k=0$).} \label{fig:Flex_ReachVolume}
%\vspace{-0.5cm}
%\end{figure}

\section{Conclusion}\label{sec:Conclusion}

The contributions of this paper are the presented modeling and analysis techniques for the quantitative assessment and visualization of operational flexibility in electric power systems.

%These techniques allow in a first phase the modeling and definition of operational flexibility of individual power system units by building up on our previous work on the Power Nodes modeling framework~\cite{PowerNodes:2010, PowerNodes:2012} and combining it with the valuable work of others, notably~\cite{Makarov:2009}. 
%In a second phase, the analysis and visualization of the operational flexibility of individual power system units is presented for some illustrative examples. The approaches are, however, also applicable for more complex, larger-scale power system setups.
%In a third phase, the aggregation of operational flexibility from several, different individual power system units is explained and illustrated. This allows notably the analysis of the overall flexibility properties of unit pools, in which different power system units are aggregated and work together to achieve a common control objective. Also, the calculation of the \emph{remaining} operational flexibility in a power system after having subtracted the \emph{needed} from the originally \emph{available} operational flexibility was shown for an illustrative case.

We envision that these techniques will become useful tools for system operators, allowing the aggregation of the available (often too) plentiful power system state information into intuitive visual charts, i.e.~3D images of \emph{available} and \emph{needed} operational flexibility cubes, and straight-forward flexibility quantification, i.e.~the flexibility metrics~$\{\rho,\pi,\epsilon\}$, for the current system state as well as for predicted future system states.

This would notably allow the real-time analysis of the overall flexibility properties of unit pools, in which different power system units are aggregated and work together to achieve a common control objective (e.g.~frequency and power balance regulation) but also the calculation of the \emph{remaining} operational flexibility in a power system after having subtracted the \emph{needed} flexibility for mitigating a disturbance (e.g.~forecast error) from the originally \emph{available} operational flexibility.

\bibliographystyle{IEEEtran}	% (uses file "plain.bst")
%\bibliography{IEEEabrv,literature}		% expects file "myrefs.bib"%
\bibliography{PSCC}		% expects file "myrefs.bib"%

\end{document}